\documentclass[leqno,12pt]{amsart} 
\usepackage{amsmath,amsfonts,amssymb,amsthm}
\usepackage{dsfont,mathrsfs,enumerate}
\usepackage{curves,epic,eepic}

\theoremstyle{plain} 
\newtheorem{thm}{\indent\sc Theorem}[section]
\newtheorem{lemma}[thm]{\indent\sc Lemma}
\newtheorem{cor}[thm]{\indent\sc Corollary}
\newtheorem{prop}[thm]{\indent\sc Proposition}

\theoremstyle{definition} 
\newtheorem{defi}[thm]{\indent\sc Definition}
\newtheorem{rem}[thm]{\indent\sc Remark}
\newtheorem{ex}[thm]{\indent\sc Example}
\newtheorem{con}[thm]{\indent\sc Construction}

\def\Z{\mathds Z}
\def\Q{\mathds Q}
\def\R{\mathds R}
\def\C{\mathds C}
\def\phi{\varphi} 
\def\id{{\rm id}}
\def\A{\mathds A}
\def\P{\mathds P}
\def\O{\mathcal O}

\DeclareMathOperator{\Mor}{Mor}

\DeclareMathOperator{\Spec}{Spec}
\DeclareMathOperator{\PGL}{PGL}
\DeclareMathOperator{\Desc}{Desc}
\DeclareMathOperator{\Pic}{Pic}
\DeclareMathOperator{\NE}{NE}
\DeclareMathOperator{\Nef}{Nef}
\DeclareMathOperator{\conv}{conv}

\begin{document}

\title[Toric varieties associated with Weyl chambers]
{The functor of toric varieties associated with Weyl chambers 
and Losev-Manin moduli spaces}

\author[Victor Batyrev]{Victor Batyrev}
\address{Mathematisches Institut, Universit\"at T\"ubingen,\newline
Auf der Morgenstelle 10, 72076 T\"ubingen, Germany}
\email{batyrev@everest.mathematik.uni-tuebingen.de}

\author[Mark Blume]{Mark Blume}
\thanks{
The second author was supported by DFG-Schwerpunkt 1388 Darstellungstheorie.
\newline\newline}
\address{Mathematisches Institut, Universit\"at M\"unster,\newline
Einsteinstrasse 62, 48149 M\"unster, Germany}
\email{mark.blume@uni-muenster.de}

\maketitle

\begin{abstract}
A root system $R$ of rank $n$ defines an $n$-dimensional
smooth projective toric variety $X(R)$ associated with 
its fan of Weyl chambers. We give a simple description of 
the functor of $X(R)$ in terms of the root system $R$ 
and apply this result in the case of root systems of type 
$A$ to give a new proof of the fact that the toric variety 
$X(A_n)$ is the fine moduli space $\overline{L}_{n+1}$ of 
stable $(n+1)$-pointed chains of projective lines 
investigated by Losev and Manin. 
\end{abstract}

\section*{Introduction} 

Let $R\subset E$ be a root system of rank $n$ in an $n$-dimensional
Euclidean space $E$ and let $M(R)\subset E$ be its root lattice.
The toric variety $X(R)$ corresponding to the root system $R$ is the
smooth projective toric variety associated with the fan of Weyl
chambers $\Sigma(R)$ in the dual space $E^*$ with respect to the
lattice $N(R)\subset E^*$ dual to the root lattice $M(R)$. 
It was shown by Klyachko in \cite{Kl85} (see also \cite{Kl95}) that if $G$ 
is a semisimple algebraic group corresponding to $R$ and $B$ is 
a Borel subgroup in $G$ then the toric variety $X(R)$ can be 
characterised as the closure of a general orbit of a maximal torus 
$T \subset G$ acting on the flag variety $G/B$.
The natural representation of the Weyl group $W(R)$ on the cohomology 
of $X(R)$ has been studied by Procesi \cite{Pr90}, 
Dolgachev-Lunts \cite{DL94}, and Stembridge \cite{St94}.

The present paper is inspired by a paper of Losev and Manin \cite{LM00},
in which fine moduli spaces $\overline{L}_n$ of stable $n$-pointed 
chains of projective lines were constructed and it was observed that 
the Losev-Manin moduli space $\overline{L}_n$ is the toric variety
associated with the polytope called the $(n-1)$-dimensional 
{\sl permutohedron} studied by Kapranov \cite[(4.3)]{Ka93}. 
These toric varieties form the $A_n$-family of the toric varieties
associated with root systems: the Losev-Manin moduli space 
$\overline{L}_{n+1}$ coincides with the toric variety $X(A_n)$.
Moreover, in \cite{LM00} is was shown that the homology groups of 
$\overline{L}_{n+1}$ $(n \geq 0)$ together with the natural action of 
the Weyl group $W(A_n)\cong S_{n+1}$ are closely related to the so called 
{\sl commutativity equations} (see also the recent paper of Shadrin and 
Zvonkine \cite{SZ09}).
We remark that these varieties are special examples of 
toric varieties obtained as equivariant blowups of $\P^n$
considered recently by Bloch and Kreimer \cite[\S 3]{BK08}.  

The Losev-Manin moduli space $\overline{L}_{n}$ is an equivariant 
compactification of a maximal torus 
$T\cong(\C^*)^{n}/\C^*\subset\PGL(n,\C)$, where the torus $T$ can 
be identified as the moduli space of $n$ points in 
$\P^1\setminus\{0,\infty\}$ up to automorphisms of $\P^1$ fixing 
$0$ and $\infty$. The boundary components of $\overline{L}_{n}$ 
parametrise certain types of $n$-pointed reducible rational curves.
There are some similarities and relations between the Losev-Manin 
moduli spaces and the well-known Grothendieck-Knudsen moduli spaces.
The Losev-Manin moduli spaces $\overline{L}_n$ parametrise
isomorphism classes of chains of projective lines with two 
poles and $n$ marked points that may coincide, whereas the 
Grothendieck-Knudsen moduli spaces $\overline{M}_{0,n+2}$ 
parametrise isomorphism classes of trees of projective lines 
with $n+2$ marked points that may not coincide.
They are related by surjective birational morphisms 
$\overline{M}_{0,n+2}\to\overline{L}_n$ dependent on the choice 
of two different elements $i,j\in\{1,\ldots,n+2\}$.
Both form a particular case of moduli spaces of weighted 
pointed stable curves as introduced by Hassett \cite{Ha03}. 

Our main objective was to generalise the result of Losev and Manin 
to other root systems $R$. This problem was mentioned by Losev and 
Manin in the introduction of their paper \cite{LM00}. 
In \cite{BB11} we present results in this direction for $R$ 
a classical root system.

To investigate interpretations of the toric varieties $X(R)$ 
associated with root systems as moduli spaces, it is natural first 
to investigate their functors of points. The functor of toric varieties 
in general was described by Mumford in \cite[Ch.\ I]{AMRT}; a different 
description was proposed by Cox \cite{Co95} for smooth toric varieties.
In the present paper we propose another description of the functor 
of the toric varieties $X(R)$ for root systems $R$, which is based on 
projection maps $X(R)\to\P^1$ and done with a view toward interpretations 
of these varieties as moduli spaces of pointed trees of projective lines.

\medskip

{\bf Outline of the paper.} 
In the first section of this paper, we derive some general results 
about the toric varieties $X(R)$ associated with arbitrary root 
systems $R$. Important are functorial properties with respect to 
maps of root systems. For example, any pair of opposite roots 
$\{\pm\alpha\}\subset R$, i.e.\ a root subsystem of type $A_1$, 
gives rise to a projection $X(R)\to X(A_1)\cong\P^1$. Morphisms 
constructed this way appear in many variants in the following.
As a main result, we give a description of the functor of the 
toric varieties $X(R)$ in terms of the root system $R$. We use 
the property of the spaces $X(R)$ that morphisms $Y\to X(R)$ are 
uniquely determined by their compositions with all the 
projection maps $X(R)\to\P^1$ given by the root subsystems 
$\{\pm\alpha\}\subset R$ of type $A_1$. Further, the relations 
between these morphisms are given by the root subsystems of type 
$A_2$ in $R$.

For the rest of the present paper, we are concerned with the toric 
varieties $X(A_n)$ associated with root systems of type $A$ and 
their interpretation as Losev-Manin moduli spaces $\overline{L}_{n+1}$.

We consider the toric varieties $X(A_n)$ in Section \ref{sec:X(A)}.
We review some results concerning the (co)homology of $X(A_n)$,
we give a basis for the homology and, in a simple way, derive the 
relations between torus invariant cycles used in \cite[Section 3]{LM00}.
Further, we comment on primitive collections and relations of the 
toric variety $X(A_n)$ and apply this to show that the 
anticanonical class of $X(A_n)$ is a semiample divisor.
This implies that $X(A_n)$ is an almost Fano variety. 
The anticanonical divisor defines a birational toric morphism 
to the Gorenstein toric Fano variety $\P_{\Delta(A_n)}$ 
corresponding to the reflexive polytope 
$\Delta(A_n)=(\textrm{convex hull of all roots of $A_n$})$.  

In Section \ref{sec:lm}, we give a new proof of the fact that 
the toric varieties $X(A_n)$ are the fine moduli spaces 
$\overline{L}_{n+1}$ of $(n+1)$-pointed chains of projective 
lines introduced by Losev and Manin.
We use the functorial properties of toric varieties associated
with root systems developed in Subsection \ref{subsec:morphisms-orbits} 
to construct the universal curve $X(A_{n+1})\to X(A_n)$ in 
Subsection \ref{subsec:univcurve}.
Our result about the functor of $X(R)$ (Subsection \ref{subsec:functor}) 
is used in the case of root systems of type $A$ in Subsection 
\ref{subsec:X(A)-L_n} to show that the functor of $X(A_n)$ is 
isomorphic to the moduli functor of $(n+1)$-pointed chains of 
projective lines.
This provides an alternative proof of the fact that this
moduli problem admits a fine moduli space $\overline{L}_{n+1}$ 
and furthermore shows that it coincides with the toric variety 
$X(A_n)$.
We will see that the data describing morphisms $Y\to X(A_n)$ 
correspond in a natural way to parameters in equations describing 
stable $(n+1)$-pointed chains of projective lines over $Y$ embedded 
in $(\P^1_Y)^{n+1}$.

\medskip
\section{Toric varieties associated with root systems and their functor}
\label{sec:X(R)}

\smallskip
\subsection{The toric variety \boldmath $X(R)$}
Let $R$ be a (reduced and crystallographic) root system in a 
Euclidean space $E$. With $R$ we associate a toric variety $X(R)$ 
(\cite{Pr90}, \cite{DL94}).

\medskip

Let $M(R)$ be the root lattice of $R$, i.e.\ the lattice in 
$E$ generated by the roots of $R$, and let $N(R)$ be the lattice
dual to $M(R)$. 
For any set of simple roots $S$, we have a cone $\sigma_S:=S^\vee
=\{v\in N(R)_\Q;\langle u,v\rangle\geq0\;\textit{for all}\;u\in S\}$ 
in the vector space $N(R)_\Q$, the (closed) Weyl chamber
corresponding to $S$.

\begin{defi} We define $\Sigma(R)$ to be the fan in the lattice 
$N(R)$ that consists of the Weyl chambers of the root system $R$ 
and all their faces. Let $X(R)$ be the toric variety associated 
with the fan $\Sigma(R)$.
\end{defi}

For $v\in N(R)_\Q$ let $\sigma_v\in\Sigma(R)$ be the cone minimal in 
$\Sigma(R)$ containing $v$. Equivalently, $\sigma_v$ is the cone 
dual to the roots in $v^\vee=\{u\in M(R)_\Q;\langle u,v\rangle
\geq 0\}$, i.e.\ $\sigma_v=(v^\vee\cap R)^\vee$.
In particular, for a general choice of $v$ this is the cone $\sigma_S$ 
dual to the set of simple roots $S\subset R$ of the set of positive 
roots $v^\vee\cap R$ defined by $v$, i.e.\ the Weyl chamber for $S$. 
The Weyl chambers cover $N(R)_\Q$, so the fan $\Sigma(R)$ is complete.
Note that each set of simple roots forms a basis of the root lattice 
$M(R)$ and $\sigma_S^\vee\cap M(R)=\langle S\rangle$ is the submonoid 
of $M(R)$ generated by $S$. $X(R)$ is covered by the open subvarieties 
$U_S:=\Spec\Z[\sigma_S^\vee\cap M(R)]=\Spec\Z[\langle S\rangle]
\cong\A^{\dim M(R)}$ for all the different sets of simple roots $S$. 

\medskip

The toric variety $X(R)$ is smooth and projective. It carries in 
a natural way the action of the Weyl group $W(R)$ of the root 
system $R$. The Weyl group permutes the sets of simple roots
and this way it acts simply transitive on the set of Weyl chambers.
The corresponding action on $X(R)$ permutes the open sets $U_S$, it
is a simply transitive action on the set of torus fixed points of 
$X(R)$.

\medskip

The root lattice $M(R)$ of $R$ is the lattice of characters of the
dense torus $T(R)$ in $X(R)$. This way, any element $u\in M(R)$ 
determines a character $x^u$ of $T(R)$, i.e.\ a rational function 
on $X(R)$.

\begin{ex}
The toric variety $X(A_1)$ is isomorphic to $\P^1$.
\end{ex}

\begin{rem}
For two root systems $R_1,R_2$ there is an isomorphism of fans
$\Sigma(R_1\times R_2)\cong\Sigma(R_1)\times\Sigma(R_2)$ and thus
an isomorphism of toric varieties $X(R_1\times R_2)\cong 
X(R_1)\times X(R_2)$.
\end{rem}

\smallskip
\subsection{Morphisms for maps of root systems and closures of 
torus orbits}\label{subsec:morphisms-orbits}
First, we show that maps between root systems coming from linear
maps of the ambient vector spaces induce toric morphisms 
of the associated toric varieties. 

\begin{prop}
Let $R,R'$ be root systems in Euclidean spaces $E,E'$. 
Then a map of vector spaces $\mu\colon E'\to E$ such that 
$\mu(R')\subset\{a\alpha;\:\alpha\in R, a\in\Z\}$ 
induces a toric morphism of the associated toric varieties 
$X(\mu)\colon X(R)\to X(R')$.
\end{prop}
\begin{proof}
The map of vector spaces $\mu\colon E'\to E$ induces a map of the
root lattices $\mu\colon M(R')\to M(R)$ because $\mu(R')\subset
\{a\alpha;\alpha\in R, a\in\Z\}$.
Let $\nu\colon N(R)\to N(R')$ be the dual map of the dual lattices.

We have to show that each cone of $\Sigma(R)$ is mapped by 
$\nu:N(R)_\Q\to N(R')_\Q$ into a cone of $\Sigma(R')$. 
Let $v\in N(R)_\Q$, we show that $\nu(\sigma_{v})\subseteq\sigma_{\nu(v)}$ 
(where as above $\sigma_{v}=(v^\vee\cap R)^\vee$ is the cone minimal 
in $\Sigma(R)$ containing $v$; in the same way the cone $\sigma_{\nu(v)}$ 
of $\Sigma(R')$ is defined). It suffices to show that 
$\mu(\nu(v)^\vee\cap R')\subseteq\langle v^\vee\cap R\rangle$. 
This is true, since $\mu(R')\subset\{a\alpha;\alpha\in R, a\in\Z\}$ 
by assumption and $\mu(\nu(v)^\vee)\subseteq v^\vee$ because 
$\langle u',\nu(v)\rangle=\langle\mu(u'),v\rangle$ for any 
$u'\in M(R')_\Q$.
\end{proof}

We have two special cases:\medskip

\noindent
(1) {\it Root subsystems induce proper surjective morphisms.}
Let $R\subset E$ be a root system and $R'\subset E'$ a root 
system in a subspace $E'\subseteq E$ such that $R'\subseteq R$.
Then $\mu\colon M(R')\to M(R)$ is injective, its dual 
$\nu\colon N(R)\to N(R')$ is surjective and we have a proper 
surjective morphism $X(R)\to X(R')$ which locally is given by 
inclusions of coordinate rings $\Z[\sigma_{\nu(v)}^\vee\cap M(R')]
\to\Z[\sigma_{v}^\vee\cap M(R)]$.

\pagebreak
\noindent
(2) {\it Projections of root systems induce closed embeddings.} 
Let $R\subset E$, $R'\subset E'$ be root systems and $\mu\colon E'\to E$ 
a homomorphism of vector spaces such that $R\subseteq\mu(R')\subset
\{a\alpha;\alpha\in R, a\in\Z\}$.
Then $\mu\colon M(R')\to M(R)$ is surjective and for $v\in N(R)_\Q$
induces a surjection $\langle\nu(v)^\vee\cap R'\rangle\to
\langle v^\vee\cap R\rangle$, the map $\nu\colon N(R)\to N(R')$ 
is injective and $\nu^{-1}(\sigma_{\nu(v)})=\sigma_{v}$ for $v\in N(R)_\Q$.
We have a closed embedding $X(R)\to X(R')$ which locally is given by 
surjective maps of coordinate rings 
$\Z[\sigma_{\nu(v)}^\vee\cap M(R')]\to\Z[\sigma_{v}^\vee\cap M(R)]$.

\begin{ex}
The first case in particular occurs if the root system $R'$ is of 
the form $R'=R\cap E'$, i.e.\ cut out by a subspace $E'\subseteq E$.
Then the morphism $X(R)\to X(R')$ is locally given by inclusions of 
coordinate rings $\Z[\sigma^\vee\!\cap M(R')]\to
\Z[\sigma^\vee\!\cap M(R)]$ for $\sigma\in\Sigma(R)$.
For example consider root subsystems $\{\pm\alpha\}\subseteq R$ 
consisting of two opposite roots, i.e.\ isomorphic to $A_1$. 
Each of these gives rise to a projection $\phi_{\{\pm\alpha\}}\colon
X(R)\to X(A_1)\cong\P^1$.
\end{ex}

\begin{ex}
For any root system $R$ the projections $X(R)\to\P^1$ for all
root subsystems $A_1\cong R'=\{\pm\alpha\}\subseteq R$ form
a morphism $X(R)\to\prod_{A_1\cong R'\subseteq R}\P^1$. 
This morphism is an instance of the second case: it 
corresponds to the projection of root systems 
$\prod_{A_1\cong R'\subseteq R}R'\to R$.
(A variant of this closed embedding has been considered 
in \cite{BJ07}.)
\end{ex}

In the second case we can describe the equations for $X(R)$
in $X(R')$.

\begin{prop}
Let $R\subset E$, $R'\subset E'$ be root systems and $\mu\colon E'\to E$ 
a homomorphism of vector spaces such that $R\subseteq\mu(R')\subset
\{a\alpha;\alpha\in R, a\in\Z\}$.
Then the image of the closed embedding $X(\mu)\colon X(R)\to X(R')$
is determined by the equations $x^{u'}=1$ for $u'\in\ker(\mu)\cap M(R')$.
Locally, the subvariety $X(R)\cap U_{S'}\subseteq U_{S'}$ for any set of 
simple roots $S'$ of $R'$ is given by the equations $\prod_ix^{\alpha_i}
=\prod_jx^{\beta_j}$ for collections of simple roots 
$\alpha_i,\beta_j\in S'$ such that 
$\sum_i\alpha_i-\sum_j\beta_j\in\ker(\mu)$.
\end{prop}
\begin{proof}
Let $v\in N(R)_\Q$ be an element in the interior of some Weyl chamber, 
let $S$ be the set of simple roots of $R$ with respect to $v$ and 
let $S'$ be the set of simple roots of $R'$ with respect to $\nu(v)$.
Then $X(\mu)^{-1}(U_{S'})=U_S$ and the inclusion $U_S\to U_{S'}$
corresponds to the surjective map of coordinate rings 
$\Z[\langle S'\rangle]\to\Z[\langle S\rangle]$ given by the surjection 
$\langle S'\rangle\to\langle S\rangle$ determined by $\mu$.
We have $\langle S\rangle\cong\langle S'\rangle/\!\!\sim$, where $\sim$ is
the equivalence relation $s_1\sim s_2\Leftrightarrow s_1-s_2\in\ker(\mu)$.
Thus $\Z[\langle S\rangle]\cong\Z[\langle S'\rangle]/I$, where $I$ is the 
ideal generated by $x^{s_1}-x^{s_2}$ for $s_1,s_2\in\langle S'\rangle$ 
such that $s_1-s_2\in\ker(\mu)$.
We can write $s_1$ (resp.\ $s_2$) as sums $s_1=\sum_i\alpha_i$ 
(resp.\ $s_2=\sum_j\beta_j$) of simple roots $\alpha_i,\beta_j\in S'$.
\end{proof}

\begin{ex}
In the case of the embedding $X(R)\to\prod_{A_1\cong R'\subseteq R}\P^1$,
the equations come from the linear relations between the roots of the 
root system $R$. This will be discussed in detail in Subsection
\ref{subsec:functor}. 
\end{ex}

Next, we consider closures of torus orbits in $X(R)$. In general, 
such orbit closures are again toric varieties and are in bijection
with the cones of the fan (see e.g.\ \cite[3.1]{Ful}). We will
see that in the case of toric varieties associated with root systems
$R$ the orbit closures are again toric varieties associated with 
certain root subsystems of $R$. 

\pagebreak
\begin{prop}\label{prop:orbitclosures}
The closure of the torus orbit $Z\subseteq X(R)$ corresponding to 
a cone $\tau\in\Sigma(R)$ is isomorphic to the toric variety $X(R')$ 
associated with the root subsystem $R'=R\cap E'\subset E'$ of $R$ 
cut out by the subspace $E'=\tau^\bot\subseteq E$ with root lattice
$M(R')=M(R)\cap E'$.

Let $S$ be a set of simple roots of $R$ such that $\tau$ is contained 
in the Weyl chamber $\sigma_S$ and put $S'=S\cap E'$. Then
$S'$ is a set of simple roots of the root system $R'$. 
The orbit closure $Z$ is covered by the open sets $Z\cap U_S$ for 
such sets of simple roots $S$, the closed subvariety 
$Z\cap U_S\subseteq U_S$ is given by the equations $x^\alpha=0$ for 
$\alpha\in S\setminus S'$.
\end{prop}
\begin{proof}
Let $S$ be a set of simple roots of $R$ such that $\tau$ is a face 
of $\sigma_S$. Then $E'=\tau^\bot$ cuts out a face of 
$\langle S\rangle_{\Q_{\geq0}}$ generated by $S'=S\cap E'$, the set
$S'$ is a set of simple roots of the root system $R'=R\cap E'$
and $M(R')=M(R)\cap E'$.

By the general theory of toric varieties (see e.g.\ \cite[3.1]{Ful}) 
the orbit closure corresponding to the cone $\tau$ is a toric variety 
covered by the affine charts $Z\cap\Spec\Z[\sigma_S^\vee\cap M(R)]
=\Spec\Z[\sigma_S^\vee\cap M(R)\cap\tau^\bot]$ for the maximal cones 
$\sigma_S\in\Sigma(R)$ such that $\tau$ is a face of $\sigma_S$. 
In the present case $\sigma_S^\vee\cap M(R)\cap\tau^\bot
=\langle S\rangle\cap\tau^\bot=\langle S'\rangle$, and so $Z$ is 
covered by the open sets $Z\cap U_S=\Spec\Z[\langle S'\rangle]$ and
isomorphic to the toric variety $X(R')$ associated with the root 
system $R'$.
The inclusion $Z\cap U_S=\Spec\Z[\langle S'\rangle]\subseteq
\Spec\Z[\langle S\rangle]=U_S$ is given by the homomorphism 
$\Z[\langle S\rangle]\to\Z[\langle S'\rangle]$, $x^u\mapsto x^u$ 
if $u\in\langle S'\rangle$ and $x^u\mapsto 0$ otherwise. Thus the 
closed subvariety $Z\cap U_S\subseteq U_S$ is determined by the 
equations $x^\alpha=0$ for $\alpha\in S\setminus S'$.
\end{proof}

Concerning the situation of the proposition we have two further remarks.

\begin{rem}\label{rem:orbitclosuresI}
The Dynkin diagram of $R'$ is the subdiagram of the Dynkin diagram of 
$R$ formed with respect to the set of simple roots $S$, that arises 
after leaving out the vertices (and adjacent edges) corresponding 
to the roots $\alpha\in S\setminus S'$. 
Usually, the root system $R'$ will be reducible and decompose as
$R'\cong\prod_iR_i$ into a number of irreducible root systems $R_i$ 
corresponding to the connected components of the Dynkin diagram of $R'$.
\end{rem}

\begin{rem}\label{rem:orbitclosuresII}
Since the fan $\Sigma(R)$ is symmetric under reflection in the 
origin, also $-\tau$ is a cone of $\Sigma(R)$ and, apart from the 
inclusion $i^+\colon X(R')\to X(R)$ of $X(R')$ as orbit closure 
corresponding to $\tau$, there is another inclusion $i^-\colon X(R')
\to X(R)$ that embeds $X(R')$ as orbit closure corresponding to $-\tau$.
Any such root subsystem $R'\subset R$ comes with a proper surjective
morphism $X(R)\to X(R')$, the two inclusions then are sections with 
respect to this morphism.

Consider the particular case of a one-dimensional cone
$\tau=\langle v\rangle_{\Q_{\geq0}}$: for $\tau,-\tau$ we have the 
two torus invariant divisors isomorphic to $X(R')$, $R'=R\cap v^\bot$, 
given as the images $i^\pm(X(R'))\subseteq X(R)$ and defined by the equations
$x^\alpha=0$ for $\alpha\in R$ such that $\langle\alpha,\pm v\rangle>0$.
\end{rem}

\smallskip
\subsection{The functor of \boldmath $X(R)$}\label{subsec:functor}
We will give a description of the functor of the toric variety 
$X(R)$ in terms of the root system $R$.
This is done via the proper surjective toric morphisms $X(R)\to\P^1$ 
for root subsystems isomorphic to $A_1$ forming the closed embedding 
$X(R)\to\prod_{A_1\cong R'\subseteq R}\P^1$ and by the use 
of the functor of $\P^1$. 

\begin{rem} We recall the well known description of the functor 
of $\P^1$.
For any scheme $Y$ a morphism $Y\to\P^1$ is uniquely determined 
(with respect to chosen coordinates on $\P^1$) by the data consisting
of a line bundle $\mathscr L$ on $Y$ together with two sections 
$t_-,t_+$ that generate $\mathscr L$ up to isomorphisms of line bundles 
with two sections (or equivalently by a line bundle $\mathscr L$ 
together with a surjective homomorphism $\O_Y^{\oplus 2}\to\mathscr L$
up to isomorphism). The functor that associates to a scheme $Y$
such data and to a morphism $Y'\to Y$ the map given by pull-back of 
line bundles and sections then is isomorphic to the functor 
$\Mor(\:\cdot\:,\P^1)$. An isomorphism is determined by the universal 
data consisting of the twisting sheaf together with homogeneous 
coordinates $(\O_{\P^1}(1),x_0,x_1)$  on $\P^1$ corresponding 
to $\id_{\P^1}\in\Mor(\P^1,\P^1)$.

Another equivalent formulation is as follows: take as data
open sets $U_-,U_+\subseteq Y$ and regular functions $f_-\in\O_Y(U_-)$,
$f_+\in\O_Y(U_+)$ such that $Y=U_-\cup U_+$, $f_-f_+=1$ on
$U_-\cap U_+$ and $\{f_-\neq 0\}=U_-\cap U_+=\{f_+\neq 0\}$.
Such data $(U_-,U_+,f_-,f_+)$ corresponds to a line bundle with 
two generating sections $(\mathscr L,t_-,t_+)$ uniquely determined
up to isomorphism by $U_-=\{t_+\neq 0\}$, $U_+=\{t_-\neq 0\}$ and 
$f_-=t_-/t_+$ on $U_-$, $f_+=t_+/t_-$ on $U_+$.
\end{rem}

By Subsection \ref{subsec:morphisms-orbits}, root subsystems of a 
root system $R$ isomorphic to $A_1$ define toric morphisms 
$X(R)\to\P^1$. These morphisms can also be described in terms of 
the preceding remark.
 
\begin{ex}
Any root subsystem of $R$ isomorphic to $A_1$, i.e.\ an unordered pair 
of opposite roots $\{\pm\alpha\}$ in $R$, defines a morphism $X(R)\to\P^1$. 
The data $\{(U_\alpha,f_\alpha),(U_{-\alpha},f_{-\alpha})\}$ are defined 
in terms of the rational functions $x^\alpha,x^{-\alpha}$ associated 
with the roots $\alpha,-\alpha$: let $U_\alpha$ be the open subset of 
$X(R)$ where $x^\alpha$ is regular, $f_\alpha:=x^\alpha|_{U_\alpha}$ 
and the same for $-\alpha$. 
The rational functions $x^\alpha,x^{-\alpha}$ have no common 
zeros or poles because any half-space in $M(R)$ contains at least 
one of the roots $\alpha,-\alpha$. 
\end{ex}

We will denote these morphisms by 
\[
\phi_{\{\pm\alpha\}}\colon\:X(R)\:\to\:\P^1_{\{\pm\alpha\}}
\]
We consider $\P^1_{\{\pm\alpha\}}$ as a copy of $\P^1$ with chosen 
homogeneous coordinates $z_{\alpha},z_{-\alpha}$. 
The toric morphism $\phi_{\{\pm\alpha\}}$ corresponds to a map of 
lattices $\Z u_\alpha\to M(R)$, $u_\alpha\mapsto\alpha$ and hence 
to a homomorphism of algebras $\Z[y_\alpha^\pm]\to\Z[M(R)]$, 
$y_\alpha\mapsto x^\alpha$, where $y_\alpha=z_{\alpha}/z_{-\alpha}$, 
i.e.\ the pull-back of the rational function 
$y_{\alpha}=z_{\alpha}/z_{-\alpha}$ on $\P^1$ via $\phi_{\{\pm\alpha\}}$ 
is the rational function $x^\alpha$ on $X(R)$. 

\medskip

Also by the preceding subsection, the collection of these morphisms
defines a closed embedding 
\[\textstyle\phi\colon\;X(R)\;\to\;\prod_{A_1\cong R'\subseteq R}
\P^1_{R'}\,.\]
We choose a set of positive roots $R^+$ of $R$. This toric morphism 
corresponds to a surjective map of lattices 
$\mu\colon\bigoplus_{\alpha\in R^+}\Z u_\alpha\to M(R)$, 
$u_\alpha\mapsto\alpha$ or of algebras 
$\bigotimes_{\alpha\in R^+}\Z[y_\alpha^{\pm}]\to\Z[M(R)]$,
$y_\alpha\mapsto x^\alpha$.
The equations describing $X(R)$ in $\prod_{A_1\cong R'\subseteq R}
\P^1_{R'}$ come from elements in the kernel of $\mu$, and such an element
$u=\sum_il_iu_{\alpha_i}\in\ker(\mu)$ corresponds to a linear relation
$\sum_il_i\alpha_i=0$ among the positive roots of the root system $R$. 
For any such element we have an equation $\prod_iy_{\alpha_i}^{l_i}=1$
or equivalently a homogeneous equation 
$\prod_iz_{\alpha_i}^{l_i}=\prod_iz_{-\alpha_i}^{l_i}$.

\begin{ex}
Consider the toric variety $X(A_2)$ associated with the root 
system $A_2=\{\pm\alpha,\pm\beta,\pm(\gamma=\alpha+\beta)\}$ 
and its embedding 
$X(A_2)\to\P^1_{\{\pm\alpha\}}\times\P^1_{\{\pm\beta\}}
\times\P^1_{\{\pm\gamma\}}$.
There is a one-dimensional space of linear relations generated
by the relation $\alpha+\beta=\gamma$, so $X(A_2)\subset
\P^1_{\{\pm\alpha\}}\times\P^1_{\{\pm\beta\}}\times
\P^1_{\{\pm\gamma\}}$ is determined by the homogeneous equation 
$z_{\alpha}z_{\beta}z_{-\gamma}=z_{-\alpha}z_{-\beta}z_{\gamma}$.
\end{ex}

In general, for any root subsystem in $R$ isomorphic to $A_2$, 
there is a linear relation of the form $\alpha+\beta=\gamma$. 
We show that these generate the space of all linear relations.  

\begin{prop}
Let $R$ be a root system. Then the space of linear relations 
between the positive roots of $R$ is generated by the relations 
$\alpha+\beta=\gamma$ for root subsystems 
$\{\pm\alpha,\pm\beta,\pm(\gamma=\alpha+\beta)\}$
of $R$ isomorphic to $A_2$.
\end{prop}
\begin{proof}
We show that the kernel of the map of lattices 
$\mu\colon\bigoplus_{\alpha\in R^+}\Z u_\alpha\to M(R)$
is generated by elements of the form $u_\alpha+u_\beta-u_\gamma$.

Since the simple roots $\alpha_i$ form a basis of the lattice $M(R)$, 
the lattice $\ker(\mu)$ is generated by elements of the form 
$u_\beta-\sum_il_iu_{\alpha_i}$, where $\beta$ is a positive 
root and thus a linear combination of simple roots 
$\beta=\sum_il_i\alpha_i$ with $l_i\in\Z_{\geq0}$.

The statement now follows from the fact that starting with the set 
of simple roots $S_0:=S$ one obtains all positive roots by successively 
adding roots that are sums of\linebreak two roots already obtained, i.e.\
$S_{i+1}=S_i\cup\{\gamma\in R;\:\textit{$\gamma=\alpha+\beta$ 
for some $\alpha,\beta\in S_i$}\}$\linebreak  
(see \cite[Ch.4, \S 1.6, Prop.19]{Bou} or \cite[\S 21.3]{FH}).
\end{proof}

\begin{cor}\label{cor:equations-X(R)}
The image of the closed embedding $\phi\colon\;X(R)\;\to\;
\prod_{A_1\cong R'\subseteq R}\P^1_{R'}$
is determined by the homogeneous equations 
$z_{\alpha}z_{\beta}z_{-\gamma}=z_{-\alpha}z_{-\beta}z_{\gamma}$ 
for root subsystems
$\{\pm\alpha,\pm\beta,\pm(\gamma=\alpha+\beta)\}$ of $R$ 
isomorphic to $A_2$.
\end{cor}

With a view to the closed embedding $\phi\colon\;X(R)\;\to\;
\prod_{A_1\cong R'\subseteq R}\P^1_{R'}$, we formulate a description 
of the functor of $X(R)$ by characterising a morphism $Y\to X(R)$ 
in terms of the family of morphisms $Y\to\P^1_{\{\pm\alpha\}}$ 
for all root subsystems of $R$ isomorphic to $A_1$ that satisfy 
compatibility conditions coming from the root subsystems of $R$ 
isomorphic to $A_2$.

\begin{defi} Let $R$ be a root system.
We define a contravariant functor\linebreak $F_R\colon(\textrm{schemes})
\to(\textrm{sets})$ that associates to a scheme $Y$ the following 
data, called $R$-data:
a family $(U_\alpha,f_\alpha)_{\alpha\in R}$ consisting of open 
sets $U_\alpha\subseteq Y$ and regular functions 
$f_\alpha\in\O_Y(U_\alpha)$ that satisfy the conditions,
\begin{enumerate}[\rm(i)]
\item $\textit{for all}\;\alpha\in R\colon\;Y=U_\alpha\cup U_{-\alpha}$, 
$\{f_\alpha\neq 0\}=U_\alpha\cap U_{-\alpha}\;\textit{and}\;
f_\alpha f_{-\alpha}=1\;\textit{on}\;U_\alpha\cap U_{-\alpha}$\,,
\item $\textit{for all}\;\alpha,\beta,\gamma\in R\colon\;\textit{if}\; 
\gamma=\alpha+\beta,\;\textit{then}\;U_{\alpha}\cap U_{\beta}
\subseteq U_\gamma\;\textit{and}\;f_{\alpha}f_{\beta}=f_\gamma
\;\textit{on}\;U_{\alpha}\cap U_{\beta}$\,,
\end{enumerate}
or, equivalently, a family $(\mathscr L_{\{\pm\alpha\}},\{t_\alpha,
t_{-\alpha}\})_{\{\pm\alpha\}\subseteq R}$ of line bundles with two 
generating sections that satisfy\\[0.3em]
\hspace*{0.3em} (ii)'
$\textit{for all}\;\alpha,\beta,\gamma\in R:\;\textit{if}\;\,
\gamma=\alpha+\beta,\;\textit{then}\;\,
t_{\alpha}t_{\beta}t_{-\gamma}=
t_{-\alpha}t_{-\beta}t_\gamma$\,,\\[0.3em] 
up to isomorphism of line bundles with a pair of sections.
To a morphism $h\colon Y'\to Y$ we associate the map 
$F_R(h)\colon F_R(Y)\to F_R(Y')$ given by pull-back of open 
sets and functions or line bundles with sections.
\end{defi}

\begin{ex}
$R$-data over a field $K$ can be written as a collection 
$((t_{\alpha}\!:\!t_{-\alpha}))_{\{\pm\alpha\}\subseteq R}$
of ratios of elements of $K$ (such that for any $\alpha$ not both 
$t_\alpha,t_{-\alpha}$ are zero) that satisfy the equations 
$t_{\alpha}t_{\beta}t_{-\gamma}=t_{-\alpha}t_{-\beta}t_\gamma$ 
for $\gamma=\alpha+\beta$.
\end{ex}

\begin{rem}\label{rem:univdataI}
On $X(R)$ we have the following $R$-data, called the universal
$R$-data, coming from the morphisms $\phi_{\{\pm\alpha\}}
\colon X(R)\to\P^1_{\{\pm\alpha\}}$, i.e.\ for $\alpha\in R$ 
consider the rational function $x^\alpha$, define $U_\alpha$ 
as the open set where $x^\alpha$ is regular and put 
$f_\alpha:=x^\alpha|_{U_\alpha}$.
\end{rem}

\begin{thm}
The toric variety $X(R)$ associated with the root system $R$ 
together with the universal $R$-data represents the functor $F_R$.
\end{thm}
\begin{proof}
We show that there is an isomorphism of functors 
$\Mor(\:\cdot\:,X(R))\cong F_R$ such that the identity in 
$\Mor(X(R),X(R))$ corresponds to the universal $R$-data on $X(R)$
denoted by $D_0\in F_R(X(R))$.

By $\Phi(Y)\colon\Mor(Y,X(R))\to F_R(Y)$, $h\mapsto h^*(D_0)$, 
we have defined a morphism of functors $\Phi\colon
\Mor(\:\cdot\:,X(R))\to F_R$.

On the other hand, for $R$-data $D$ on $Y$ we have a morphism
$\phi_D\colon Y\to X(R)\subseteq\prod_{\{\pm\alpha\}}\P^1_{\{\pm\alpha\}}$, 
where $X(R)$ is considered as a closed subvariety of 
$\prod_{\{\pm\alpha\}}\P^1_{\{\pm\alpha\}}$ via the embedding 
$\phi\colon\:X(R)\:\to\:\prod_{\{\pm\alpha\}}\P^1_{\{\pm\alpha\}}$.
In particular it is $\phi_{D_0}=\id_{X(R)}$.
The maps $\Psi(Y)\colon F_R(Y)\to\Mor(Y,X(R))$, $D\mapsto\phi_D$ 
form a morphism of functors $\Psi\colon F_R\to\Mor(\:\cdot\:,$ $X(R))$, 
because for any morphism of schemes $h\colon Y'\to Y$ it is 
$\phi_{h^*(D)}=\phi_D\circ h$. This is true since the maps
$\phi_{D,\{\pm\alpha\}}\colon Y\to\P^1_{\{\pm\alpha\}}$ associated 
with the part of $R$-data for pairs of roots $\{\pm\alpha\}$ 
satisfy $\phi_{h^*(D),\{\pm\alpha\}}=\phi_{D,\{\pm\alpha\}}\circ h$.

Thus, we have two morphisms of functors $\Phi$ and $\Psi$; these
are inverse to each other.
\end{proof}

\begin{rem}\label{rem:univdataII}
The universal $R$-data on $X(R)$ gives rise to $R$-data 
$((t_{\alpha}:t_{-\alpha}))_{\{\pm\alpha\}\subseteq R}$
over points of $X(R)$ having the following properties:\\
$\bullet$ Over the affine chart $U_S$ for a set of simple roots 
$S$ of $R$, we have $(t_{\alpha}:t_{-\alpha})\neq (1:0)$ if 
$\alpha\in\langle S\rangle$, and over the torus fixed point of $U_S$, 
we have $(t_{\alpha}:t_{-\alpha})=(0:1)$ for $\alpha\in\langle S\rangle$.\\
$\bullet$ Over the torus invariant divisor corresponding to a one-dimensional
cone generated by $v$, we have $(t_{\alpha}:t_{-\alpha})=(0:1)$ if
$\langle\alpha,v\rangle>0$ (cf.\ Remark \ref{rem:orbitclosuresII}).
\end{rem}

\medskip
\section{Toric varieties associated with root systems of type $A$}
\label{sec:X(A)}

\smallskip
\subsection{Toric varieties \boldmath $X(A_n)$}\label{subsec:X(A_n)}
Consider an $(n+1)$-dimensional Euclidean vector space with basis 
$u_1,\ldots,u_{n+1}$. The root system $A_n$ in the $n$-dimensional
subspace $E=\{\sum_ia_iu_i;\linebreak\sum_ia_i=0\}$ consists of the 
$n(n+1)$ roots 
\[u_i-u_j\quad\textit{for}\;\;i,j\in\{1,\ldots,n+1\},\;i\neq j.\]
The lattice $N(A_n)\cong\Z^n$ dual to the root lattice $M(A_n)\cong\Z^n$ 
has a generating system $v_1,\ldots,v_{n+1}$ with one relation
$\sum_iv_i=0$, it is a quotient of the lattice dual to 
$\bigoplus_{i=1}^{n+1}\Z u_i$ with basis $v_1,\ldots,v_{n+1}$ 
dual to $u_1,\ldots,u_{n+1}$.

\medskip

The sets of simple roots of the root system $A_n$ are of the form 
\[S=\{u_{i_1}-u_{i_2},u_{i_2}-u_{i_3},\ldots,u_{i_{n}}-u_{i_{n+1}}\}\] 
for some ordering $i_1,\ldots,i_{n+1}$ of the set $\{1,\ldots,{n+1}\}$. 
The maximal cone $\sigma_S=S^\vee$ of $\Sigma(A_n)$, i.e.\ the Weyl 
chamber corresponding to $S$, consists of those elements 
$\sum_ia_iv_i\in N(A_n)$ that satisfy 
$a_{i_1}\geq a_{i_2},\ldots,a_{i_n}\geq a_{i_{n+1}}$ or 
equivalently of non-negative linear combinations of 
\[v_{i_1},\;v_{i_1}+v_{i_2},\;\ldots\,,\;v_{i_1}+\cdots+v_{i_{n}}.\] 

\medskip

So, the fan $\Sigma(A_n)$ can be described as follows (and coincides 
with that of \cite[(2.5), (2.6)]{LM00}): 
there are $2^{n+1}-2$ one-dimensional cones of the fan $\Sigma(A_n)$, 
and these are generated by the elements $v_A=\sum_{i\in A}v_i\in N(A_n)$ 
for $A\in\mathcal A$ where 
\[\mathcal A=\{A\:|\:\emptyset\neq A\subsetneq\{1,\ldots,n+1\}\}.\]
A family $(v_{A^{(i)}})_{i=1,\ldots,k}$ corresponding to a collection 
of pairwise different sets $A^{(1)},\ldots,A^{(k)}$ $\in\mathcal A$ 
generates a $k$-dimensional cone of $\Sigma(A_n)$ whenever these sets 
can be ordered such that $A^{(i_1)}\subset A^{(i_2)}\subset\cdots
\subset A^{(i_k)}$.

\medskip

The fan $\Sigma(A_n)$ defines an $n$-dimensional smooth projective
toric variety $X(A_n)$. It is covered by the $(n+1)!$ open subvarieties 
$U_S=\Spec\Z[\sigma_S^\vee\cap M(A_n)]$ for the sets of simple roots 
$S$ corresponding to strict orderings of the set $\{1,\ldots,n+1\}$.
If $S=\{u_{i_1}-u_{i_2},u_{i_2}-u_{i_3},\ldots,u_{i_{n}}-u_{i_{n+1}}\}$,
then $\Z[\sigma_S^\vee\cap M(A_n)]=\Z[x_{i_1}/x_{i_2},\ldots,
x_{i_{n}}/x_{i_{n+1}}]$.
The Weyl group of the root system $A_n$ is the symmetric group $S_{n+1}$;
it acts on $\Sigma(A_n)$ and on $X(A_n)$ by permuting the Weyl chambers 
and the open sets $U_S$.

\medskip

By results of the last subsection (Proposition \ref{prop:orbitclosures} 
and Remark \ref{rem:orbitclosuresI}), the closures of torus orbits in 
$X(A_n)$ are isomorphic to products $\prod_i X(A_{n_i})$. In particular, 
the torus invariant divisor corresponding to the cone generated by 
$v_{i_1}+\cdots+v_{i_k}$ is of the form $X(A_{n-k})\times X(A_{k-1})$.

\pagebreak
\begin{figure}[h]
\begin{picture}(150,30)(-5,0)

\put(-45,0){
\put(80,15){\vector(1,0){14}}\put(80,15){\vector(-1,0){14}}
\put(96,15){\makebox(0,0)[l]{$v_1$}}
\put(64,15){\makebox(0,0)[r]{$v_2$}}
\drawline(80,14)(80,16)
}

\put(2,23){\makebox(0,0)[l]{\large$\Sigma(A_1)$}}
\put(148,23){\makebox(0,0)[r]{\large$\Sigma(A_2)$}}

\put(25,0){
\put(80,15){\vector(1,0){14}}\put(80,15){\vector(-1,0){14}}
\put(80,15){\vector(1,2){7}}\put(80,15){\vector(1,-2){7}}
\put(80,15){\vector(-1,2){7}}\put(80,15){\vector(-1,-2){7}}
\put(90,29){\makebox(0,0)[l]{$v_1+v_2$}}
\put(90,1){\makebox(0,0)[l]{$v_1+v_3$}}
\put(70,29){\makebox(0,0)[r]{$v_2$}}
\put(70,1){\makebox(0,0)[r]{$v_3$}}
\put(97,15){\makebox(0,0)[l]{$v_1$}}
\put(63,15){\makebox(0,0)[r]{$v_2+v_3$}}
}
\end{picture}
\vspace{-1cm}
\caption{}
\end{figure}
\vspace{-3mm}

\begin{ex}
The fans $\Sigma(A_1)$ and $\Sigma(A_2)$ are as Figure 1. 
\end{ex}

There are other descriptions of the toric variety $X(A_n)$:

\medskip

{\it $X(A_n)$ as blow-up of $\P^n$.}
$X(A_n)$ can be constructed by a sequence of toric blow-ups starting 
with $\P^n$ by first blowing up the $n+1$ torus fixed points, then 
blowing up the strict transform of the lines joining two of these 
points, then the strict transform of the planes through any three
of these points and so on (see \cite[Ch.3]{Pr90}, \cite[(4.3.13)]{Ka93}, 
\cite[(5.1)]{DL94}).

\medskip

{\it $X(A_n)$ as toric variety associated with the 
lattice polytope called permutohedron.}
The $n$-dimensional permutohedron is defined as the convex 
hull in $\Q^{n+1}$ of the $S_{n+1}$-orbit of the point
$(1,2,\ldots,n+1)$ where the symmetric group $S_{n+1}$ 
acts by permuting coordinates.
One can show that the fan for this  polytope is the fan  
$\Sigma(A_n)$. The toric variety $X(A_n)$ also appears as 
``permutohedral space $\Pi^n\,$'' 
in \cite[(4.3.10) through (4.3.13)]{Ka93}.

\smallskip
\subsection{(Co)homology of \boldmath $X(A_n)$}\label{subsec:cohomX(A_n)}
We know that the integral cohomology is torsion-free and confined 
to the even degrees. Standard methods from toric geometry 
(see e.g.\ \cite[(10.8)]{Dan}) furthermore imply the
description of the cohomology ring of the toric variety $X(A_n)$ 
over the complex numbers (see also \cite[(2.7)]{LM00}) as
\[H^*(X(A_n),\Z)\;\cong\;\Z[\,l_A;A\in\mathcal A\,]/(R_1+R_2),\]
where $R_1$ is the ideal generated by the elements
$r_{i,j}=\sum_{i\in A}l_A-\sum_{j\in A}l_A$ for $i,j\in\{1,\ldots,n+1\}$, 
$i\neq j$, and $R_2$ the ideal generated by the elements 
$r_{A,A'}=l_Al_{A'}$ for $A,A'\in\mathcal A$ such that 
$A\not\subseteq A'$, $A'\not\subseteq A$ (these correspond to
the primitive collections of the fan $\Sigma(A_n)$, see next 
subsection).

\medskip

The Betti numbers and Poincar\'e polynomials of the varieties $X(A_n)$ 
over the complex numbers are calculated in \cite[(2.3)]{LM00}; see 
also \cite[Section 6]{St94}, \cite[Section 4]{DL94} for a description 
in terms of the Eulerian numbers and see \cite{Kl85}, \cite{Kl95},
\cite{St94} for the general case of toric varieties associated 
with root systems. In particular,  we know that the rank of 
$H^*(X(A_n),\Z)$ is $(n+1)!$, i.e.\ the number of maximal cones
of $\Sigma(A_n)$.

\medskip

The $\Z$-module $\Z[\,l_A:A\in\mathcal A\,]/(R_1+R_2)$ is generated
by the classes of square-free monomials (see \cite[(10.7.1)]{Dan}).
We can restrict to monomials each of which has only factors 
corresponding to one-dimensional faces of one maximal cone. 
Such a monomial $\prod_{i=1}^ml_{A^{(i)}}$ corresponds to an 
$m$-dimensional face of the respective maximal cone and, on the 
other hand, to a collection $A^{(1)}\subsetneq\cdots\subsetneq A^{(m)}$ 
of elements of $\mathcal A$.
We denote by $G$ the $\Z$-submodule of $\Z[\,l_A;A\in\mathcal A\,]$
generated by these monomials (called ``good monomials'' in 
\cite[(2.8)]{LM00}). We have the canonical isomorphism of $\Z$-modules 
$G/U\cong\Z[\,l_A;A\in\mathcal A\,]/(R_1+R_2)$ where 
$U=(R_1+R_2)\cap G$.
The module $G/U$ can be identified with the homology module 
$H_*(X(A_n),\Z)$ (cf.\ \cite[(2.9.2)]{LM00}).

In \cite[(2.8.2)]{LM00} the following generators of the module of 
relations $U$ are given. For a collection $A^{(1)}\subsetneq\cdots
\subsetneq A^{(m)}$ of elements of $\mathcal A$ and 
$k\in\{1,\ldots,m+1\}$, $i,j\in A^{(k)}\setminus A^{(k-1)}$ (put
$A^{(0)}=\emptyset$, $A^{(m+1)}=\{1,\ldots,n+1\}$), $i\neq j$, let 
\[\textstyle
r_{i,j}((A^{(h)})_h,k)
=\Big(\sum_{\tiny\hspace{-1mm}
\begin{array}{l}i\!\in\!A\\[-1mm]j\!\not\in\!A\\\end{array}
\hspace{-1.2mm}}l_A-
\sum_{\tiny\hspace{-1mm}
\begin{array}{l}j\!\in\!A\\[-1mm]i\!\not\in\!A\\\end{array}
\hspace{-1.2mm}}l_A\Big)
\prod_{h=1}^ml_{A^{(h)}}
\]
where the sums run over sets $A\in\mathcal A$ such that 
$A^{(k-1)}\subsetneq A\subsetneq A^{(k)}$.

The maximal cones of the fan $\Sigma(A_n)$ correspond to collections 
$A^{(1)}\subsetneq\cdots\subsetneq A^{(n)}$ of elements of 
$\mathcal A$ and these correspond to permutations $\sigma\in S_{n+1}$ 
via $\{\sigma(1),\ldots,\sigma(k)\}=A^{(k)}$ for $k=1,\ldots,n$.
The descent set of a permutation $\sigma\in S_{n+1}$ is the set 
\[\Desc(\sigma)=\{k\in\{1,\ldots,n\};\:\sigma(k)>\sigma(k+1)\}\,.\]
For any $\sigma\in S_{n+1}$ we define a monomial in $G$ by
\[\textstyle
l^\sigma=\prod_{k\not\in\Desc(\sigma)}l_{\{\sigma(1),\ldots,\sigma(k)\}}\,.
\]

\enlargethispage{3mm}
\begin{prop}
The classes of the monomials $l^\sigma$ for $\sigma\in S_{n+1}$ form 
a basis of the homology module $G/U=H_*(X(A_n),\Z)$. 
The module of relations $U$ is generated by the elements 
$r_{i,j}((A^{(h)})_h,k)$ for collections $A^{(1)}\subsetneq
\ldots\subsetneq A^{(m)}$ of elements of $\mathcal A$ and 
$k\in\{1,\ldots,m+1\}$, $i,j\in A^{(k)}\setminus A^{(k-1)}$ 
{\rm(}put $A^{(0)}=\emptyset$, $A^{(m+1)}=\{1,\ldots,n+1\}${\rm)}, 
$i\neq j$. 
\end{prop}
\begin{proof}
We have $(n+1)!$ distinct monomials $l^\sigma$, and this number coincides 
with the rank of $G/U$. Thus it remains to show that every monomial in 
$G$ via the relations $r_{i,j}((A^{(h)})_h,k)$ is equivalent to a 
linear combination of the monomials $l^\sigma$.

For a monomial $\prod_{k=1}^ml_{A^{(k)}}$ corresponding to a collection
$A^{(1)}\subsetneq\cdots\subsetneq A^{(m)}$ we define the number
$d(\prod_{k=1}^ml_{A^{(k)}}):=|\{k\in\{1,\ldots,m\};\:
\min P_k>\max P_{k+1}\}|\in\Z_{\geq 0}$ in terms of the associated 
partition $P_1=A^{(1)}$, $P_2=A^{(2)}\setminus A^{(1)}$, $\ldots\;$, 
$P_m=A^{(m)}\setminus A^{(m-1)}$, $P_{m+1}=\{1,\ldots,n+1\}\setminus 
A^{(m)}$ of the set $\{1,\ldots,n+1\}$. The monomials $y\in G$ satisfying 
$d(y)=0$ are exactly the monomials of the form $l^\sigma$.
On the other hand we have an ordering $\prec$ of the monomials of $G$:
take the partition $(P_i)_{i=1,\ldots,m+1}$ associated with a monomial
and consider the sequence (corresponding to a permutation of the set
$\{1,\ldots,n+1\}$) that arises by taking first the elements of $P_{m+1}$ 
then those of $P_m$ and so on and by ordering the elements of each 
$P_i$ according to their size, and on these sequences we take the 
lexicographic order.

We show that every monomial in $G$ modulo $U$ is equivalent to a 
linear combination of the monomials $l^\sigma$, $\sigma\in S_{n+1}$, 
by showing that every monomial $y\in G$ with $d(y)>0$ modulo a 
relation is equivalent to a linear combination of monomials $y'$ 
with $y\prec y'$.
In fact, let $A^{(1)}\subsetneq\cdots\subsetneq A^{(m)}$ be a 
collection of elements of $\mathcal A$ with associated partition 
$(P_k)_{k=1,\ldots,m+1}$ such that the corresponding monomial
$y=\prod_{k=1}^ml_{A^{(k)}}$ satisfies $d(y)>0$.
Take $k\in\{1,\ldots,m\}$ such that $i:=\min P_k>\max P_{k+1}=:j$, 
then 
\[\textstyle
r_{i,j}((A^{(h)})_{h\neq k},k)=
\Big(\sum_{\tiny\hspace{-1mm}
\begin{array}{l}i\!\in\!A\\[-1mm]j\!\not\in\!A\\\end{array}
\hspace{-1.2mm}}l_A-
\sum_{\tiny\hspace{-1mm}
\begin{array}{l}j\!\in\!A\\[-1mm]i\!\not\in\!A\\\end{array}
\hspace{-1.2mm}}l_A\Big)
\prod_{h\neq k}l_{A^{(h)}}\,,
\]
where the sums run over sets $A\in\mathcal A$ such that $A^{(k-1)}
\subsetneq A\subsetneq A^{(k+1)}$, is a relation that contains $y$ 
as the monomial minimal with respect to $\prec$.
\end{proof}

\begin{rem}
The above set of generators for the module of relations is used in 
\cite[Section 3]{LM00} in the context of solutions to the 
commutativity equations; the respective statement \cite[(2.9)]{LM00} 
is proven in a different way in that paper.
Our statement involves a basis of homology which, in the general case
of toric varieties associated with root systems, is given in 
\cite{Kl85}, \cite{Kl95} (see \cite[Rem.\ 5.4]{BB11} for more 
explanations).
\end{rem}

\subsection{Primitive collections and the morphism \boldmath 
$X(A_n)\to\P_{\Delta(A_n)}$}
It was observed in \cite{Ba91} that any $n$-dimensional smooth 
projective toric variety $X$ corresponding to a fan $\Sigma$ 
can be described by the set of primitive collections among the
generators of the one-dimensional cones of $\Sigma$ together with 
the corresponding primitive relations.
A primitive collection  of the fan $\Sigma$ is a set of generators of 
one-dimensional cones that does not generate a cone of $\Sigma$, but all 
of its proper subsets generate a cone of $\Sigma$.

\begin{thm} {\rm(\cite[Thm.\ 2.15]{Ba91}, see also 
\cite[Thm.\ 1.4, Prop.\ 1.10]{CR08}).}
The Mori cone $\NE(X)\subset A_1(X)\otimes_\Z\R$ of effective $1$-cycles
is generated by the primitive relations.
Equivalently, the nef-cone $\Nef(X)\subset\Pic(X)\otimes_\Z\R$ {\rm(}which
coincides with the closure of the ample cone{\rm)} is given by line bundles 
corresponding to piecewise linear functions $\phi$ satisfying 
$\phi(w_1)+\cdots+\phi(w_k)\geq\phi(w_1+\cdots+w_k)$ for all primitive 
collections $\{w_1,\ldots,w_k\}$ of the fan $\Sigma$.
\end{thm}

In the case $\Sigma=\Sigma(A_n)$, the primitive collections 
consist of two elements $v_A,v_{A'}$ (again put 
$v_A=\sum_{i\in A}v_i$) corresponding to non comparable elements 
of $\mathcal A$, i.e.\ elements $A,A'\in\mathcal A$ such that 
$A\not\subseteq A'$, $A'\not\subseteq A$. The corresponding 
primitive relation has one of the following four forms
(put $I=\{1,\ldots,n+1\}$):

(1) $v_A+v_{A'}=0$, {\it if $A\cup A'=I$ and $A\cap A'=\emptyset$};

(2) $v_A+v_{A'}=v_{A\cup A'}$, {\it if $A\cup A'\neq I$ and
$A\cap A'=\emptyset$};

(3) $v_A+v_{A'}=v_{A\cap A'}$, {\it if $A\cup A'=I$ and
$A\cap A'\neq\emptyset$};

(4) $v_A+v_{A'}=v_{A\cap A'}+v_{A\cup A'}$, {\it if
$A\cup A'\neq I$ and $A\cap A'\neq\emptyset$.}

\begin{cor}\label{cor:(semi)ample}
Let $D=\sum_{A\in\mathcal A}a_AD_A$ be a torus-invariant divisor 
in $X(A_n)$, where $D_A$ is the torus-invariant prime divisor 
corresponding to $v_A$. We put $a_I=a_{\emptyset}=0$. 
Then $D$ is ample {\rm(}resp.\ semiample, or equivalently nef{\rm)} 
if and only if 
\[a_A+a_{A'}>a_{A\cap A'}+a_{A\cup A'}\]
{\rm(}resp.\ $\geq${\rm)} for all primitive collections $\{v_A,v_{A'}\}$. 
\end{cor} 

\pagebreak
\begin{cor}
The anticanonical class of $X(A_n)$ is semiample, or equivalently nef;
$X(A_n)$ is an almost Fano variety.
\end{cor} 
\begin{proof}
By corollary \ref{cor:(semi)ample} $-K_{X(A_n)}$ is semiample 
(and ample if there is no primitive relation of type (4), i.e.\ 
if $n\leq 2$). That $X(A_n)$ is an almost Fano variety means that 
$-K_{X(A_n)}$ is semiample and big, where the second property follows 
from the fact that the corresponding polytope is full-dimensional.
\end{proof}

Being semiample and big, the anticanonical class defines a toric 
morphism $X(A_n)\to\P(H^0(\O_{X(A_n)}(-K_{X(A_n)})))$, which is birational 
onto its image, but not necessarily an embedding. 

\begin{prop}
The anticanonical class defines a birational toric morphism 
$X(A_n)\to\P_{\Delta(A_n)}$, where $\P_{\Delta(A_n)}$ is 
the Gorenstein toric Fano variety associated with the reflexive polytope 
\[
\Delta(A_n)\;=\;\{m\in M_\Q;\:\left<m,v_A\right>\geq-1\;\textit{for}\;
A\in\mathcal A\}\;=\;\conv\{\alpha_{ij};\:i,j\in I, i\neq j\}
\]
and $\alpha_{ij}=u_{i}-u_{j}$ are the roots of the root system $A_n$.
\end{prop}
\begin{proof}
In general, the image of the morphism given by a semiample divisor 
$\sum_ia_iD_i$ (sum over the torus invariant prime divisors $D_i$)
is the polarised toric variety corresponding to the polytope 
$\{m\in M_\Q;\:\left<m,v_i\right>\geq-a_i\;\textit{for all}\;i\}$
($v_i\in N$ being the lattice generators of the one dimensional 
cones corresponding to $D_i$).
In the present case, we have the anticanonical divisor
$-K_{X(A_n)}=\sum_{A\in\mathcal A}D_A$ on $X(A_n)$ giving rise
to the polytope $\Delta(A_n)=\{m\in M_\Q;\:\left<m,v_A\right>\geq-1
\;\textit{for}\;A\in\mathcal A\}$.
This polytope $\Delta(A_n)$ coincides with the convex hull of
the roots of $A_n$. Clearly, every root $\alpha_{ij}$ satisfies 
$\left<m,v_A\right>\geq-1$ for all $A\in\mathcal A$.
On the other hand, a given element
$m=(m_1,\ldots,m_{n+1})\in\Delta(A_n)$ we can write as 
$m=\sum_{i\in B,j\in C}a_{ij}\alpha_{ij}$ with $a_{ij}\geq 0$ 
and $B=\{i;\:m_i>0\}$, $C=\{i;\:m_i<0\}$, and the condition 
$\left<m,v_C\right>\geq-1$ gives $\sum a_{ij}\leq 1$, so 
$m\in\conv\{\alpha_{ij}\:|\:i\neq j\}$.
The polytope $\Delta(A_n)$ is reflexive ($0$ being its only inner lattice 
point, the other lattice points in $\Delta(A_n)$ are the roots), 
equivalently $\P_{\Delta(A_n)}$ is a Gorenstein toric Fano variety.
\end{proof}

\begin{rem}
The fan $\Sigma_{\Delta(A_n)}$ determined by $\Delta(A_n)$ 
consists of the  cones
\[
\sigma_{B_1,B_2}=\left<v_A;\:B_1\subseteq A\subseteq B_2\right>
\]
for $B_1,B_2\in\mathcal A$, $B_1\subseteq B_2$. 
The cone $\sigma_{B_1,B_2}$ has dimension $|B_2\setminus B_1|+1$ and
is generated by $2^{|B_2\setminus B_1|}$ elements $v_A$.
We see that $\P_{\Delta(A_n)}$ is singular for $n\geq 3$, the
singular locus consisting of the torus orbits of codimension 
at least $3$.
The morphism $X(A_n)\to\P_{\Delta(A_n)}$ is a crepant resolution
(an MPCP-desingularisation in the sense of \cite{Ba94}). It is given 
by subdividing each $d$-dimensional cone $\sigma_{B_1,B_2}$
into $(d-1)!$ $d$-dimensional cones generated by 
$v_{B_1},v_{B_1\cup\{a_1\}},\;\ldots\;,v_{B_1\cup\{a_1,\ldots,a_{d-1}\}}
\!=\!v_{B_2}$ corresponding to permutations of the set $B_2\setminus B_1$.
A set $\{v_{A^{(1)}},\ldots,v_{A^{(k)}}\}$ is contained in a single 
cone of $\Sigma_{\Delta(A_n)}$, if $\bigcap_iA^{(i)}\neq\emptyset$ 
and $\bigcup_iA^{(i)}\neq I$. In particular, the primitive collections 
of types (1), (2) and (3) survive.
\end{rem}

\pagebreak
\section{Losev-Manin moduli spaces}
\label{sec:lm}

\smallskip
\subsection{The moduli functor of \boldmath $A_n$-curves}
\label{subsec:lmmodulispaces}
We begin with a presentation of the Losev-Manin moduli spaces 
$\overline{L}_n$ introduced in \cite{LM00} and the corresponding moduli 
functor of stable $n$-pointed chains of projective lines (here called 
$A_{n-1}$-curves).

\medskip

Consider a complex projective line $\P^1_\C$ with two distinct 
closed points $s_-,s_+\in\P^1_\C$ called poles. The moduli space 
of $n$ distinguishable points in $\P^1_\C\setminus\{s_-,s_+\}$ is 
a torus $(\C^*)^n/\C^*\cong(\C^*)^{n-1}$, here we divide out 
the automorphism group $\C^*$ of the projective line with 
two poles $(\P^1_\C,s_-,s_+)$.
This space has a natural compactification $\overline{L}_n$ such 
that its boundary parametrises isomorphism classes of certain 
types of reducible $n$-pointed curves.

\begin{defi} 
A {\sl chain of projective lines} of length $m$ over an algebraically 
closed field $K$ is a projective curve $C=C_1\cup\cdots\cup C_m$ over 
$K$ such that each irreducible component $C_j$ of $C$ is a projective 
line with poles $p^-_j,p^+_j$ and these components intersect as follows: 
different components $C_i$ and $C_j$ intersect only if $|i-j|=1$ 
and in this case $C_j,C_{j+1}$ intersect transversally at the 
single point $p^+_j=p^-_{j+1}$. 
For $p^-_1\in C_1$ (resp.\ $p^+_m\in C_m$) we write $s_-$ (resp.\ $s_+$).
Two chains of projective lines $(C,s_-,s_+)$ and $(C',s_-',s_+')$ 
are called isomorphic if there is an isomorphism $\phi\colon C\to C'$ 
such that $\phi(s_-)=s_-'$ and $\phi(s_+)=s_+'$.
\end{defi}

\begin{defi} 
An {\sl $n$-pointed chain of projective lines} 
$(C,s_-,s_+,s_1,\ldots,s_n)$ is a chain of projective 
lines together with closed points $s_i\in C$ different 
from $s_-,s_+$ and the intersection points of components (see Figure 2).
Two $n$-pointed chains of projective lines
$(C,s_-,s_+,s_1,\ldots,s_n)$ and $(C',s_-',s_+',s'_1,\ldots,s'_n)$
are called isomorphic if there is an isomorphism $\phi\colon
(C,s_-,s_+)\to(C',s_-',s_+')$ of the underlying chains of projective 
lines such that $\phi(s_j)=s'_j$ for all $j\in\{1,\ldots,n\}$.
An $n$-pointed chain of projective lines $(C,s_-,s_+,s_1,\ldots,s_n)$ 
is called {\sl stable} if each component of $C$ contains at least one 
of the points $s_j$.
An {\sl $A_n$-curve} over an algebraically closed field $K$ is a stable 
$(n+1)$-pointed chain of projective lines over $K$.

\vspace{-1mm}
\begin{figure}[h]
\begin{picture}(150,22)(-5,0)
\put(15,0){\line(2,1){40}}
\put(75,0){\line(-2,1){40}}
\put(55,0){\line(2,1){40}}
\put(55,0){\line(2,1){40}}
\put(115,0){\line(-2,1){40}}
\put(95,0){\line(2,1){40}}

\filltype{white}
\put(21,3){\circle*{1.2}}
\put(18,4){\makebox(0,0)[r]{\small$s_-\!=\!p^-_1$}}
\put(45,15){\circle*{1.2}}
\put(45,19){\makebox(0,0)[b]{\small$p^+_1\!=\!p^-_2$}}
\put(65,5){\circle*{1.2}}
\put(65,1){\makebox(0,0)[t]{\small$p^+_2\!=\!p^-_3$}}
\put(85,15){\circle*{1.2}}
\put(85,19){\makebox(0,0)[b]{\small$p^+_3\!=\!p^-_4$}}
\put(105,5){\circle*{1.2}}
\put(105,1){\makebox(0,0)[t]{\small$p^+_4\!=\!p^-_5$}}
\put(129,17){\circle*{1.2}}
\put(132,16){\makebox(0,0)[l]{\small$p^+_5\!=\!s_+$}}

\filltype{black}
\put(31,8){\circle*{1.2}}
\put(30,9){\makebox(0,0)[b]{\small$s_{i_1}$}}
\put(41,13){\circle*{1.2}}
\put(41,11){\makebox(0,0)[t]{\small$s_{i_2}$}}
\put(53,11){\circle*{1.2}}
\put(53,9){\makebox(0,0)[t]{\small$s_{i_3}$}}
\put(71,8){\circle*{1.2}}
\put(70,9){\makebox(0,0)[b]{\small$s_{i_4}$}}
\put(79,12){\circle*{1.2}}
\put(82,11){\makebox(0,0)[t]{\small$s_{i_5}$}}
\put(99,8){\circle*{1.2}}
\put(101,9){\makebox(0,0)[b]{\small$s_{i_6}$}}
\put(115,10){\circle*{1.2}}
\put(118,9){\makebox(0,0)[t]{\small$s_{i_7}$}}
\end{picture}
\vspace{-1mm}
\caption{}
\end{figure}
\end{defi}
\vspace{-6mm}

\begin{defi}\label{defi:combtype}
An $A_n$-curve of length $m$ decomposes into irreducible components 
$C=C_1\cup\cdots\cup C_m$ with $s_-\in C_1$, $s_+\in C_m$. We will call
the resulting decomposition  
\[\{1,\ldots,n+1\}=P_1\sqcup\cdots\sqcup P_m\]
such that $i\in P_k$ if and only if $s_i\in C_k$ the 
{\sl combinatorial type} of the $A_n$-curve $C$.
We will also write this in the form $s_{i_1}\ldots
s_{i_l}|s_{i_{l+1}}\ldots|\ldots$ with the sections for the different 
sets $P_k$ separated by the symbol ``$\,|\,$'' starting on the left with 
$P_1$.
\end{defi}

\begin{defi} 
Let $Y$ be a scheme. An {\sl $A_n$-curve over $Y$} is a collection 
$(\pi\colon C\to Y,s_-,s_+,s_1,\ldots,s_{n+1})$, 
where $C$ is a scheme, $\pi$ is a proper flat morphism of schemes
and $s_-,s_+,s_1,$ $\ldots,s_{n+1}\colon Y\to C$ are 
sections such that for any geometric point $y$ of $Y$ the 
collection $(C_y,(s_-)_y,(s_+)_y,(s_1)_y,$ $\ldots,(s_{n+1})_y)$ 
is an $A_n$-curve over $y$.
An isomorphism of $A_n$-curves over $Y$ is an isomorphism of 
$Y$-schemes that is compatible with the sections.
We define the {\sl moduli functor of $A_n$-curves} as the 
contravariant functor
\[
\begin{array}{crcl}
\overline{L}_{n+1}:&(\textrm{schemes})&
\to&(\textrm{sets})\\
&Y&\mapsto&\left\{\textrm{$A_n$-curves over $Y$}\right\}\,/\sim
\end{array}
\]
that associates to a scheme $Y$ the set of isomorphism classes 
of $A_n$-curves over $Y$ and to a morphism of schemes the map 
obtained by pulling back $A_n$-curves. 
\end{defi}

It is shown in \cite[(2.2)]{LM00} that the functor 
$\overline{L}_{n+1}$ is represented by a smooth projective variety,
denoted by the same symbol $\overline{L}_{n+1}$, that 
is, we have a fine moduli space of $A_n$-curves.
Further it is argued in \cite[(2.6.3)]{LM00} that 
$\overline{L}_{n+1}$ is isomorphic to the toric variety 
$X(A_n)$ associated with the root system $A_n$.

\medskip

In \cite[(2.1)]{LM00} an inductive construction of $\overline{L}_{n+1}$ 
together with the universal curve $C_{n+1}\to\overline{L}_{n+1}$
is given using arguments of \cite{Kn83}.
As in the case of the similar moduli spaces $\overline{M}_{0,n}$, 
the universal curve $C_{n+1}$ over $\overline{L}_{n+1}$ is isomorphic 
to the next moduli space $\overline{L}_{n+2}$. 

\begin{ex}\label{ex:univcurve}\ \\
(1) $C_1\to\overline{L}_1$ is isomorphic to $\P^1\to pt$ with 
$0,\infty$ and $1$-section. This reflects the fact that 
any $A_0$-curve over a scheme $Y$ is isomorphic to the trivial 
projective bundle $\P^1_Y$ with $0,\infty$ and $1$-section.\\
(2) $\overline{L}_2$ is isomorphic to $\P^1$. The fibres of the 
universal curve $C_2$ over $\P^1\setminus\{(0:1),(1:0)\}$ are 
pointed $(\P^1,s_-,s_+)$, and over $(0:1)$ and $(1:0)$ the fibres 
are pointed chains consisting of two components (see Figure 3).

\bigskip
\begin{figure}[h]
\setlength{\unitlength}{1.45mm}
\begin{picture}(150,40)(20,0)
\put(31,27){\makebox(0,0)[l]{$C_2$}}
\put(31,5){\makebox(0,0)[l]{$\overline{L}_2$}}

\put(43,5){\line(1,0){64}}
\put(51,4){\line(0,1){2}}\put(51,4){\makebox(0,0)[t]{\Small$(0:1)$}}
\put(99,4){\line(0,1){2}}\put(99,4){\makebox(0,0)[t]{\Small$(1:0)$}}

\put(51,22){\line(-1,4){2}}\put(51,22){\line(1,-4){2}}
\put(51,32){\line(1,4){2}}\put(51,32){\line(-1,-4){2}}
\put(99,22){\line(1,4){2}}\put(99,22){\line(-1,-4){2}}
\put(99,32){\line(-1,4){2}}\put(99,32){\line(1,-4){2}}
\dottedline{0.8}(43,18)(107,18)\dottedline{0.8}(43,36)(107,36)
\put(59,37){\makebox(0,0)[b]{\Small$s_+$}}
\put(59,17){\makebox(0,0)[t]{\Small$s_-$}}
\put(59,30){\makebox(0,0)[t]{\Small$s_1$}}
\put(59,22){\makebox(0,0)[t]{\Small$s_2$}}
\curve[8](43,22, 51,22, 59,23, 67,25, 75,27, 83,29, 91,31, 
99,32, 107,32)
\curve[8](43,32, 51,32, 59,31, 67,29, 75,27, 83,25, 91,23, 
99,22, 107,22)
\end{picture}
\vspace{-1cm}
\caption{}
\setlength{\unitlength}{1mm}
\end{figure}
\end{ex}

\pagebreak
\subsection{The universal curve}\label{subsec:univcurve}
We construct the universal $A_n$-curve $X(A_{n+1})\to X(A_{n})$ 
with its sections, which later will be seen to coincide with the 
universal curve over the Losev-Manin moduli space $\overline{L}_{n+1}$,
in terms of the functorial properties of toric varieties associated
with root systems developed in subsection \ref{subsec:morphisms-orbits}.

\begin{con}[The universal $A_n$-curve]\label{con:univcurve}
Consider the root subsystem $A_n\subset A_{n+1}$
consisting of the roots $\beta_{ij}=u_i-u_j$ for 
$i,j\in\{1,\ldots,n+1\}$. The inclusion of root systems 
$A_n\subset A_{n+1}$ determines a proper surjective morphism 
$X(A_{n+1})\to X(A_n)$.

There are the $n+1$ additional pairs of opposite roots $\pm\alpha_i$, 
$\alpha_i=u_i-u_{n+2}$ for $i\in\{1,\ldots,n+1\}$.
For each of the root subsystems $\{\pm\alpha_i\}\cong A_1$ in $A_{n+1}$
not contained in $A_n$ we have a section $s_i\colon X(A_n)\to X(A_{n+1})$
corresponding to the projection of the root system $A_{n+1}$
onto the root subsystem $A_n$ with kernel generated by $\alpha_i$.
The image of $s_i$ can be described by the equation $x^{\alpha_i}=1$.
On the other hand, the inclusion $\{\pm\alpha_i\}\subset A_{n+1}$ 
defines a projection $X(A_{n+1})\to X(A_1)\cong\P^1_{\{\pm\alpha_i\}}$.  
If we choose coordinates of $\P^1_{\{\pm\alpha_i\}}$ as in Subsection 
\ref{subsec:functor}, then the image of the section $s_i$ is given by 
the preimage of the point $(1:1)\in\P^1_{\{\pm\alpha_i\}}$.

Further, we have two sections $s_\pm\colon X(A_n)\to X(A_{n+1})$
which are inclusions of $X(A_n)$ into $X(A_{n+1})$ as torus 
invariant divisors (cf.\ Remark \ref{rem:orbitclosuresII}). 
The image of $s_-$ (resp.\ $s_+$) corresponds to the ray of the fan 
$\Sigma(A_{n+1})$ generated by $-v_{n+2}$ (resp.\ $v_{n+2}$) and is 
given by the equations $x^{\alpha_i}=0$ (resp.\ $x^{-\alpha_i}=0$). 
\end{con}

\begin{prop} \
The collection
$(X(A_{n+1})\!\to\!X(A_n),s_-,s_+,s_1,\ldots,s_{n+1})$ 
is an $A_n$-curve over $X(A_n)$.
\end{prop}
\begin{proof}
The morphism $X(A_{n+1})\to X(A_n)$ is proper. 
To show flatness, consider the covering by affine toric charts: 
for any set of simple roots $S$ of the root system $A_n$ 
there are $n+2$ affine spaces $U_{S_j}\cong\A^{n+1}$ lying over
$U_S\cong\A^{n}$ corresponding to $n+2$ sets of simple roots 
$S_j$ of the root system $A_{n+1}$ such that the submonoid 
$\langle S_j\rangle$ of the root lattice $M(A_{n+1})$ generated 
by $S_j$ contains $S$. For example, if 
$S=\{u_1-u_2,\ldots,u_n-u_{n+1}\}$ we have 
$S_j=\{\ldots,u_{j-1}-u_j,u_j-u_{n+2},u_{n+2}-u_{j+1},
u_{j+1}-u_{j+2},\ldots\}$ for $j=0,\ldots,n+1$ 
(in particular $S_0=\{u_{n+2}-u_1,u_1-u_2,\ldots\}$ and
$S_{n+1}=\{\ldots,u_{n}-u_{n+1},u_{n+1}-u_{n+2}\}$).
The maps $U_{S_j}\to U_S$ are flat morphisms $\A^{n+1}\to\A^n$. 
Indeed, they are given by the identity on a polynomial ring in $n-1$ 
variables tensored with a map of the form 
$\Z[z]\to\Z[x,y]$, $z\mapsto xy$ or $z\mapsto y$.

By what is shown below, the fibres are $A_n$-curves:
Remark \ref{rem:emb-univcurve} describes the universal 
$A_n$-curve in terms of equations given by $A_n$-data, and
Proposition \ref{prop:emb-A_n-curve} shows that such equations 
define an $A_n$-curve once it is known that any $A_n$-data over 
a field arises as in Proposition \ref{prop:emb-A_n-curve} from an 
$A_n$-curve. This is shown in Lemma \ref{le:A_n-data--A_n-curve}.
\end{proof}

\begin{defi}
We call the object $(X(A_{n+1})\!\to\!X(A_n),s_-,s_+,s_1,\ldots,s_{n+1})$ 
of Construction \ref{con:univcurve} the universal $A_n$-curve over 
$X(A_n)$.
\end{defi}

\begin{ex} 
The universal curve $C_2$ over $\overline{L}_2$ is pictured as
Figure 3 and later seen to coincide with $X(A_2)$ over $X(A_1)$. 
Here, in Figure 4, we draw the corresponding inclusion of root systems 
(left) and the map of fans $\Sigma(A_2)\to\Sigma(A_1)$ (right).

\vspace{-2mm}
\begin{figure}[h]
\begin{picture}(150,60)(-5,0)
\put(3,40){\makebox(0,0)[l]{\large$A_2$}}
\put(3,5){\makebox(0,0)[l]{\large$A_1$}}
\put(35,40){\vector(1,0){14}}\put(35,40){\vector(-1,0){14}}
\put(35,40){\vector(1,2){7}}\put(35,40){\vector(-1,2){7}}
\put(35,40){\vector(1,-2){7}}\put(35,40){\vector(-1,-2){7}}
\dottedline{1}(14,44)(56,44)\dottedline{1}(14,36)(56,36)
\dottedline{1}(14,44)(14,36)\dottedline{1}(56,44)(56,36)
\put(49,37){\makebox(0,0)[b]{\small$u_1\!-\!u_2$}}
\put(21,37){\makebox(0,0)[b]{\small$u_2\!-\!u_1$}}
\put(46,54){\makebox(0,0)[l]{\small$\alpha_1\!=\!u_1\!-\!u_3$}}
\put(24,54){\makebox(0,0)[r]{\small$\alpha_2\!=\!u_2\!-\!u_3$}}
\put(46,26){\makebox(0,0)[l]{\small$-\alpha_2\!=\!u_3\!-\!u_2$}}
\put(24,26){\makebox(0,0)[r]{\small$-\alpha_1\!=\!u_3\!-\!u_1$}}
\put(26,33){\vector(0,1){2}}\put(44,33){\vector(0,1){2}}
\dottedline{1}(26,8)(26,35)\dottedline{1}(44,8)(44,35)
\put(35,5){\vector(1,0){14}}\put(35,5){\vector(-1,0){14}}
\put(49,1){\makebox(0,0)[b]{\small$u_1\!-\!u_2$}}
\put(21,1){\makebox(0,0)[b]{\small$u_2\!-\!u_1$}}
\drawline(35,4)(35,6)
\put(147,40){\makebox(0,0)[r]{\large$\Sigma(A_2)$}}
\put(147,5){\makebox(0,0)[r]{\large$\Sigma(A_1)$}}
\put(110,40){\vector(0,1){14}}\put(110,40){\vector(0,-1){14}}
\put(110,40){\vector(2,1){14}}\put(110,40){\vector(-2,1){14}}
\put(110,40){\vector(2,-1){14}}\put(110,40){\vector(-2,-1){14}}
\put(110,24){\makebox(0,0)[t]{\small$v_3$}}
\put(126,50){\makebox(0,0)[l]{\small$v_1$}}
\put(94,30){\makebox(0,0)[r]{\small$v_2+v_3$}}
\put(110,56){\makebox(0,0)[b]{\small$v_1+v_2$}}
\put(126,30){\makebox(0,0)[l]{\small$v_1+v_3$}}
\put(94,50){\makebox(0,0)[r]{\small$v_2$}}
\dottedline{1}(100,28)(100,12)\put(100,12){\vector(0,-1){2}}
\dottedline{1}(120,28)(120,12)\put(120,12){\vector(0,-1){2}}
\put(110,5){\vector(1,0){14}}\put(110,5){\vector(-1,0){14}}
\put(126,5){\makebox(0,0)[l]{\small$v_1$}}
\put(94,5){\makebox(0,0)[r]{\small$v_2$}}
\drawline(110,4)(110,6)
\end{picture}
\vspace{-2mm}
\caption{}
\end{figure}
\vspace{-3mm}
\end{ex}

We apply the embedding into a product of projective lines 
\[\textstyle
X(R)\to\prod_{A_1\cong R'\subseteq R}\P^1_{R'}=:P(R)
\]
considered in Subsections \ref{subsec:morphisms-orbits} and 
\ref{subsec:functor} to the situation of Construction 
\ref{con:univcurve}.

\begin{rem}\label{rem:emb-univcurve}
Consider $X(A_{n+1})$ (resp.\ $X(A_{n})$) as embedded into 
the product of projective lines $P(A_{n+1})$ (resp.\ $P(A_n)$).
Then the morphism $X(A_{n+1})\to X(A_n)$ is induced by the projection 
onto the subproduct $P(A_{n+1})\to P(A_n)$.
The subvarieties $X(A_{n+1})\subseteq P(A_{n+1})$ (resp.\ 
$X(A_{n})\subseteq P(A_n)$) are determined by the homogeneous equations 
$z_{\alpha}z_{\beta}z_{-\gamma}=z_{-\alpha}z_{-\beta}z_{\gamma}$
for roots $\alpha,\beta,\gamma$ such that $\alpha+\beta=\gamma$,
i.e.\ root subsystems of type $A_2$ in $A_{n+1}$ (resp.\ $A_n$). 
If we first consider the product $P(A_{n+1})$ and the equations coming 
from root subsystems of type $A_2$ in $A_n$, this gives 
\[\textstyle
P(A_{n+1}/A_n)_{X(A_n)}=\big(\prod_{A_1\cong R\subseteq 
A_{n+1}\setminus A_n}\P^1_R\big)_{X(A_n)}
=\big(\prod_{i=1}^{n+1}\P^1_{\{\pm\alpha_i\}}\big)_{X(A_n)}\,.
\]
Therein, $X(A_{n+1})$ is the closed subscheme given by 
the homogeneous equations
\begin{equation}\label{eq:univ-A_n-curve}
t_{\beta_{ij}}z_{\alpha_j}z_{-\alpha_i}
=t_{-\beta_{ij}}z_{-\alpha_j}z_{\alpha_i},
\quad i,j\in\{1,\ldots,n+1\}, i\neq j,\;\beta_{ij}=\alpha_i-\alpha_j
\end{equation}
where $t_{\pm\beta_{ij}}$ are the homogeneous coordinates of 
$\P^1_{\{\pm\beta_{ij}\}}$ (consider $X(A_{n})$ as embedded into 
$\prod_R\P^1_R$) or equivalently the two generating 
sections of the line bundle $\mathscr L_{\{\pm\beta_{ij}\}}$ being 
part of the universal $A_n$-data on $X(A_{n})$.
The sections $s_i\colon X(A_n)\to X(A_{n+1})$ for $i\in\{1,\ldots,n+1\}$ 
are given by the additional equations $z_{\alpha_i}=z_{-\alpha_i}$. 
The sections $s_-$ (resp.\ $s_+$) are given by $z_{\alpha_1}=\;\cdots\;=
z_{\alpha_{n+1}}=0$ (resp.\ $z_{-\alpha_1}=\;\cdots\;=z_{-\alpha_{n+1}}=0$).
\end{rem}

\begin{ex}
The universal $A_1$-curve $X(A_2)\subset(\P^1_{\{\pm\alpha_1\}}
\times\P^1_{\{\pm\alpha_2\}})_{X(A_1)}$ 
over $X(A_1)$ is given by the homogeneous equation 
\[
t_{\beta_{12}}z_{\alpha_2}z_{-\alpha_1}
=t_{-\beta_{12}}z_{-\alpha_2}z_{\alpha_1}
\]
where $\beta_{12}=u_1-u_2$ and
$(\mathscr L_{\{\pm\beta_{12}\}},\{t_{\beta_{12}},t_{-\beta_{12}}\})$ 
is the universal $A_1$-data on $X(A_1)\cong\P^1$. 
We picture in Figure 5 the $A_1$-curves in $\P^1\times\P^1$ defined 
by this equation over the points given by 
$(t_{\beta_{12}}:t_{-\beta_{12}})=(0:1),(a:b),(1:0)$. In the
second case, we write the marked points in terms of the homogeneous 
coordinates $(z_{-\alpha_2}:z_{\alpha_2})$.

\vspace{-4mm}
\begin{figure}[h]
\begin{picture}(150,60)(-5,0)
\put(5,6){\line(1,0){140}}
\put(25,0){\makebox(0,0)[b]
{\small$(t_{\beta_{12}}:t_{-\beta_{12}})=(1:0)$}}
\put(25,5){\line(0,1){2}}
\dottedline{0.8}(15,15)(15,45)
\put(35,25){\line(0,1){30}}
\put(10,16){\line(2,1){30}}
\dottedline{0.8}(10,39)(40,54)
\filltype{white}
\put(15,18.5){\circle*{1.2}}
\put(13,20){\makebox(0,0)[r]{\small$s_-$}}
\put(35,51.5){\circle*{1.2}}
\put(37,50){\makebox(0,0)[l]{\small$s_+$}}
\dottedline{1.6}(25,20)(25,50)
\dottedline{1.6}(10,27.5)(40,42.5)
\filltype{black}
\put(35,40){\circle*{1.2}}
\put(37,39){\makebox(0,0)[l]{\small$s_2$}}
\put(25,23.5){\circle*{1.2}}
\put(27,22){\makebox(0,0)[l]{\small$s_1$}}
\put(75,0){\makebox(0,0)[b]
{\small$(t_{\beta_{12}}:t_{-\beta_{12}})=(a:b)$}}
\put(75,5){\line(0,1){2}}
\dottedline{0.8}(65,15)(65,45)
\dottedline{0.8}(85,25)(85,55)
\dottedline{0.8}(60,16)(86,29)
\dottedline{0.8}(60,39)(86,52)
\filltype{white}
\put(65,18.5){\circle*{1.2}}
\put(63,20){\makebox(0,0)[r]{\small$s_-$}}
\put(85,51.5){\circle*{1.2}}
\put(87,52){\makebox(0,0)[l]{\small$s_+\!=\!(0\!:\!1)$}}
\dottedline{1.6}(75,20)(75,50)
\dottedline{1.6}(60,27.5)(86,40.5)
\filltype{black}
\qbezier(60.5,14)(80,31.75)(86.3,56)
\put(78.7,36.8){\circle*{1.2}}
\put(79,39){\makebox(0,0)[r]{\small$s_2$}}
\put(87,41){\makebox(0,0)[l]{\small$s_2\!=\!(1\!:\!1)$}}
\put(87,37){\makebox(0,0)[l]{\small$s_1\!=\!(a\!:\!b)$}}
\put(75,31){\circle*{1.2}}
\put(76,30){\makebox(0,0)[l]{\small$s_1$}}
\dottedline{1.6}(75,31)(86,36.5)
\put(87,29){\makebox(0,0)[l]{\small$s_-\!=\!(1\!:\!0)$}}
\put(125,0){\makebox(0,0)[b]
{\small$(t_{\beta_{12}}:t_{-\beta_{12}})=(0:1)$}}
\put(125,5){\line(0,1){2}}
\put(115,15){\line(0,1){30}}
\dottedline{0.8}(135,25)(135,55)
\dottedline{0.8}(110,16)(140,31)
\put(110,39){\line(2,1){30}}
\filltype{white}
\put(115,18.5){\circle*{1.2}}
\put(113,20){\makebox(0,0)[r]{\small$s_-$}}
\put(135,51.5){\circle*{1.2}}
\put(137,50){\makebox(0,0)[l]{\small$s_+$}}
\dottedline{1.6}(125,20)(125,50)
\dottedline{1.6}(110,27.5)(140,42.5)
\filltype{black}
\put(115,30){\circle*{1.2}}
\put(113,31){\makebox(0,0)[r]{\small$s_2$}}
\put(125,46.5){\circle*{1.2}}
\put(124,48){\makebox(0,0)[r]{\small$s_1$}}
\end{picture}
\vspace{-4mm}
\caption{}
\end{figure}
\vspace{-4mm}
\end{ex}

By Remark \ref{rem:emb-univcurve} the universal $A_n$-curve 
over $X(A_n)$ can be embedded into a product $P(A_{n+1}/A_n)_{X(A_n)}
\cong (\P^1)^{n+1}_{X(A_n)}$ and the embedded curve is given by 
equations (\ref{eq:univ-A_n-curve}) determined by the universal 
$A_n$-data.
We show that any $A_n$-curve $C$ over a field can be embedded into
a product $(\P^1)^{n+1}$ and extract $A_n$-data such that $C$
is described by the same equations as the universal curve.

\begin{prop}\label{prop:emb-A_n-curve}
Let $(C,s_-,s_+,s_1,\ldots,s_{n+1})$ be an $A_n$-curve 
over a field. For $i\in\{1,\ldots,n+1\}$ let $z_{\alpha_i},z_{-\alpha_i}$
be a basis of $H^0(C,\O_C(s_i))$ such that $z_{\alpha_i}(s_-)=0$,
$z_{-\alpha_i}(s_+)=0$,  $z_{\alpha_i}(s_i)=z_{-\alpha_i}(s_i)\neq0$.
We will write $\P_{\{\pm\alpha_i\}}$ for $\P(H^0(C,\O_C(s_i)))$.
Then by 
\[
(t_{\beta_{ij}}:t_{-\beta_{ij}})=(z_{-\alpha_j}(s_i):z_{\alpha_j}(s_i))
\]
for $\beta_{ij}=\alpha_i-\alpha_j$, $i\neq j$, we can define $A_n$-data 
$(t_{\beta_{ij}}:t_{-\beta_{ij}})_{\{\pm\beta_{ij}\}\subseteq A_n}$,
and the morphism
\[\textstyle
C\;\to\;\prod_{i=1}^{n+1}\P^1_{\{\pm\alpha_i\}}=P(A_{n+1}/A_n)
\] 
is an isomorphism onto the closed subvariety $C'\subseteq P(A_{n+1}/A_n)$ 
determined by the homogeneous equations
\begin{equation}\label{eq:A_n-curve}
t_{\beta_{ij}}z_{\alpha_j}z_{-\alpha_i}
=t_{-\beta_{ij}}z_{-\alpha_j}z_{\alpha_i},
\quad i,j\in\{1,\ldots,n+1\}, i\neq j,\;\beta_{ij}=\alpha_i-\alpha_j\,.
\end{equation}
Furthermore, $C'$ together with the marked points $s_i'$ 
defined by the additional equations $z_{\alpha_i}=z_{-\alpha_i}$ 
and $s_-'$ {\rm(}resp.\ $s_+'${\rm)} defined by $z_{\alpha_i}=0$ for
$i=1,\ldots,n+1$ {\rm(}resp.\ $z_{-\alpha_i}=0$ for $i=1,\ldots,n+1${\rm)} 
is an $A_n$-curve and we have an isomorphism of $A_n$-curves
$(C,s_-,s_+,s_1,\ldots,$ $s_{n+1})\to(C',s_-',s_+',s_1',\ldots,s_{n+1}')$.
\end{prop}
\begin{proof}
The data $(t_{\beta_{ij}}:t_{-\beta_{ij}})$ is defined as the position 
of a marked point $s_i$ relative to another marked point $s_j$ 
of $C$. Note that the two ways to define $(t_{\beta_{ij}}:t_{-\beta_{ij}})$
are equivalent because $(z_{\alpha_i}(s_j):z_{-\alpha_i}(s_j))
=(z_{-\alpha_j}(s_i):z_{\alpha_j}(s_i))$.

The data $(t_{\beta_{ij}}:t_{-\beta_{ij}})$ relates the bases 
$z_{\alpha_i},z_{-\alpha_i}$ and $z_{\alpha_j},z_{-\alpha_j}$ 
via the equation
$t_{\beta_{ij}}z_{\alpha_j}z_{-\alpha_i}$
$=t_{-\beta_{ij}}z_{-\alpha_j}z_{\alpha_i}$.
If the corresponding marked points are not contained in the
same component of $C$, both sides of the equation are zero.
Otherwise, both are homogeneous coordinates of the same component
of $C$ and we have the formula 
$(t_{\beta_{ij}}z_{-\alpha_i}:t_{-\beta_{ij}}z_{\alpha_i})
=(z_{-\alpha_j}:z_{\alpha_j})$ which can be checked at the two poles
and the marked point $s_i$.
It follows that the morphism $C\to\P(A_{n+1}/A_n)$ maps $C$ to the 
subscheme $C'$ defined by the equations in (\ref{eq:A_n-curve}).

The collection 
$(t_{\beta_{ij}}:t_{-\beta_{ij}})_{\{\pm\beta_{ij}\}\subseteq A_n}$ 
is $A_n$-data: we have to show that the equations 
$t_{\beta_{ij}}t_{\beta_{jk}}t_{-\beta_{ik}}
=t_{-\beta_{ij}}t_{-\beta_{jk}}t_{\beta_{ik}}$
are satisfied.
If $s_i,s_j,s_k$ are not contained in the same component, then
both sides are zero. If $s_i,s_j,s_k$ are in the 
same component, then $(t_{\beta_{ik}}:t_{-\beta_{ik}})
=(z_{-\alpha_k}(s_i):z_{\alpha_k}(s_i))
=(t_{\beta_{jk}}z_{-\alpha_j}(s_i):t_{-\beta_{jk}}z_{\alpha_j}(s_i))
=(t_{\beta_{jk}}t_{\beta_{ij}}:t_{-\beta_{jk}}t_{-\beta_{ij}})$ 
making use of $(z_{-\alpha_k}:z_{\alpha_k})=
(t_{\beta_{jk}}z_{-\alpha_j}:t_{-\beta_{jk}}z_{\alpha_j})$.

We show that $C\to C'\subseteq P(A_{n+1}/A_n)$ is an isomorphism. 
The curve $C$ decomposes into irreducible components 
$C=C_{1}\cup\ldots\cup C_m$. Let $\{1,\ldots,n+1\}=P_1\sqcup
\ldots\sqcup P_m$ be the combinatorial type of $C$ (see Definition 
\ref{defi:combtype}).
For $k\in\{1,\ldots,m\}$ the morphism $C_k\to\prod_{i\in P_k}
\P^1_{\{\pm\alpha_i\}}=:P(P_k)$ is an isomorphism onto 
$C''_k\subseteq P(P_k)$ given by the equations in (\ref{eq:A_n-curve}) 
involving only coordinates $z_{\pm\alpha_i}$ for $i\in P_k$.
The equations in (\ref{eq:A_n-curve}) involving coordinates  
for roots $\alpha_i,\alpha_j$ for $i\in P_k$, $j\in P_{k'}$, 
$k\neq k'$, are of the form $z_{\alpha_j}z_{-\alpha_i}=0$ 
if $k<k'$ and $z_{-\alpha_j}z_{\alpha_i}=0$ if $k'<k$.
So, the equations in (\ref{eq:A_n-curve}) containing coordinates 
for some $i\in P_k$ define a subvariety of 
$P(A_{n+1}/A_n)=\prod_{i\in P_j,j<k}\P^1_{\{\pm\alpha_i\}}
\times P(P_k)\times\prod_{i\in P_j,k<j}\P^1_{\{\pm\alpha_i\}}$ 
consisting of the irreducible components 
$C'_k=\prod_{i\in P_{k'},k'<k}\{z_{-\alpha_i}\!=\!0\}\times C''_k
\times\prod_{i\in P_{k'},k<k'}\{z_{\alpha_i}\!=\!0\}$, 
$\prod_{i\in P_j,j<k}\P^1_{\{\pm\alpha_s\}}\times
\prod_{i\in P_j,k\leq j}\{z_{\alpha_i}=0\}$, 
$\prod_{i\in P_j,j\leq k}\{z_{-\alpha_i}=0\}\times
\prod_{i\in P_j,k<j}\P^1_{\{\pm\alpha_s\}}$.
For each $k$ we have an isomorphism $C_k\to C'_k$, and these 
form the isomorphism $C\to C'$.

The sections of $C'$ were defined such that the isomorphism 
$C\to C'$ is an isomorphism of $A_n$-curves.
\end{proof}

\begin{lemma}\label{le:A_n-data--A_n-curve}
Any $A_n$-data over a field arises as $A_n$-data extracted from
an $A_n$-curve by the method of Proposition {\rm\ref{prop:emb-A_n-curve}}.
\end{lemma}
\begin{proof} Let 
$(t_{\beta_{ij}}:t_{-\beta_{ij}})_{\{\pm\beta_{ij}\}\subseteq A_n}$ 
be $A_n$-data over a field. 

We define an ordering $\prec$ on the set $\{1,\ldots,n+1\}$. 
For $i\neq j$ define $i\prec j$ (resp.\ $i\preceq j$) if 
$(t_{\beta_{ij}}:t_{-\beta_{ij}})=(1:0)$ (resp.\ 
$(t_{\beta_{ij}}:t_{-\beta_{ij}})\neq(0:1)$) where 
$\beta_{ij}=\alpha_i-\alpha_j$. 
This ordering defines a decomposition 
$\{1,\ldots,n+1\}=P_1\sqcup\cdots\sqcup P_m$ 
into nonempty equivalence classes such that $i\prec j$ if and only if 
$i\in P_k,\,j\in P_{k'}$ for $k<k'$. 

We construct an $A_n$-curve such that the $A_n$-data extracted from 
it by the method of Proposition \ref{prop:emb-A_n-curve} is the 
given $A_n$-data.
Take a chain of projective lines $(C,s_-,s_+)$ of length $m$.
For each component $C_k$ we can choose for each $i\in P_k$ a point 
$s_i\in C_k$ different from the poles $p_k^\pm$, such that their relative 
positions are given by the data $(t_{\beta_{ij}}:t_{-\beta_{ij}})$, 
i.e.\ $s_i=(t_{\beta_{ij}}:t_{-\beta_{ij}})$ with respect to coordinates 
of $C_k$ such that $p_k^-=(1:0)$, $p_k^+=(0:1)$, $s_j=(1:1)$.
This is possible, the compatibility is assured by the conditions
of $A_n$-data. 
\end{proof}

Considering the universal $A_n$-curve, the combinatorial types of 
the geometric fibres determine a stratification of $X(A_n)$ which 
coincides with the stratification of this toric variety into torus 
orbits.

\begin{prop}\label{prop:combtypeAn}
Over the torus orbit in $X(A_n)$ corresponding to the one-dimensional 
cone generated by $v_{i_1}+\cdots+v_{i_k}$ we have the combinatorial type 
\[s_{i_{n+1}}\ldots s_{i_{k+1}}|s_{i_k}\ldots s_{i_1}\,.\]
\end{prop}
\begin{proof}
The universal $A_n$-data over each point of the closure of the orbit 
corresponding to a generator of a one-dimensional cone generated by 
$v$ has the property $(t_{\beta}:t_{-\beta})=(0:1)$ if 
$\langle\beta,v\rangle>0$ (see Remark \ref{rem:univdataII}). 
For $v=v_{i_1}+\cdots+v_{i_k}$ this in particular implies  
$(t_{\beta_{ij}}:t_{-\beta_{ij}})=(0:1)$ if $i\in\{i_1,\ldots,i_k\}$
and $j\in\{i_{k+1},\ldots,i_{n+1}\}$.
We obtain for points in this torus orbit the above combinatorial type.
\end{proof}

The proposition describes the combinatorial types over the torus 
orbits in $X(A_n)$ of codimension one. The combinatorial type over a 
lower-dimensional torus orbit is given by the partition that arises 
as common refinement of the partitions for the torus orbits of 
codimension one that contain the respective orbit in their closure.

\smallskip
\subsection{Isomorphism between the functor of \boldmath $X(A_n)$ 
and the moduli functor}
\label{subsec:X(A)-L_n}
We show that the moduli functor of $A_n$-curves is isomorphic to 
the functor of $X(A_n)$ of Subsection \ref{subsec:functor}. 
This implies that the toric variety $X(A_n)$ is a fine 
moduli space of $A_n$-curves and hence coincides with 
the Losev-Manin moduli space $\overline{L}_{n+1}$.

\medskip

To relate $A_n$-curves to $A_n$-data, we consider an embedding 
of arbitrary $A_n$-curves over a scheme $Y$ into a product 
$(\P^1)_Y^{n+1}$ that generalises the embedding in Proposition 
\ref{prop:emb-A_n-curve} to the relative situation. 
The main tools are the following contraction morphisms.
The technical details are standard, compare to similar 
constructions in \cite{Kn83}, \cite{GHP88}.

\begin{con}\label{con:contr}
Let $(\pi\colon C\to Y,s_-,s_+,s_1,\ldots,s_{n+1})$ 
be an $A_n$-curve over $Y$. Each subset of sections 
$\{s_{i_1},\ldots,s_{i_k}\}\subseteq\{s_1,\ldots,s_n\}$ defines
a line bundle $\O_C(s_{i_1}+\cdots+s_{i_k})$ on $C$ and a locally 
free sheaf $\pi_*\O_C(s_{i_1}+\cdots+s_{i_k})$ of rank $k+1$ on $Y$. 
The natural surjective homomorphism $\pi^*\pi_*\O_C(s_{i_1}+\cdots+s_{i_k})
\to\O_C(s_{i_1}+\cdots+s_{i_k})$ induces a morphism 
$p_{\{i_1,\ldots,i_k\}}\colon C\to C_{\{i_1,\ldots,i_k\}}
\subseteq\P_Y(\pi_*\mathscr\O_C(s_{i_1}+\cdots+s_{i_k}))$.

This construction commutes with base change. In particular,
for each fibre $C_y$, $y\in Y$, the morphism 
$(p_{\{i_1,\ldots,i_k\}})_y\colon C_y\to
\P_Y(\pi_*\mathscr\O_C(s_{i_1}+\cdots+s_{i_k}))_y$ arises by 
applying the above construction to the line bundle 
$\O_{C_y}(s_{i_1}(y)+\cdots+s_{i_k}(y))$ on $C_y$; it is an 
isomorphism on the components containing at least one of the
marked points $s_{i_j}(y)$ and contracts all other components.
The image of $p_{\{i_1,\ldots,i_k\}}$ with the sections
$p_{\{i_1,\ldots,i_k\}}\circ s_\pm$, $p_{\{i_1,\ldots,i_k\}}
\circ s_{i_j}$ is an $A_{k-1}$-curve.

These morphisms are functorial with respect to multiple
inclusions of sets of sections.
\end{con}

We will make use of the following particular cases:

\smallskip

(1) {\it Contraction with respect to one section onto a $\P^1$-bundle.}
For each section $s_i$ of $C$ there is a contraction morphism 
$p_i\colon C\to C_i=\P_Y(\pi_*\mathscr\O_C(s_i))$.
Since the $\P^1$-bundle $\P_Y(\pi_*\mathscr\O_C(s_i))$ has 
three disjoint sections $p_i\circ s_-,p_i\circ s_+,p_i\circ s_i\colon 
Y\to\P_Y(\pi_*\mathscr\O_C(s_i))$, there is an isomorphism 
$\P_Y(\pi_*\mathscr\O_C(s_i))\cong\P^1_Y$ that identifies 
$p_i\circ s_-,p_i\circ s_+,p_i\circ s_i$ with the 
$(1:0),(0:1),(1:1)$-section of $\P^1_Y$ (see Figure 6). 

\vspace{-3mm}
\begin{figure}[h]
\begin{picture}(150,44)(-5,0)
\put(10,22){\makebox(0,0)[l]{\large$C_y$}}
\put(30,10){\line(1,2){4}}\put(30,10){\line(-1,-2){4}}
\put(30,22){\line(-1,2){4}}\put(30,22){\line(1,-2){4}}
\put(30,34){\line(1,2){4}}\put(30,34){\line(-1,-2){4}}
\filltype{white}
\put(27,4){\circle*{1}}\put(26,4){\makebox(0,0)[r]{\small$s_-$}}
\put(33,40){\circle*{1}}\put(32,40){\makebox(0,0)[r]{\small$s_+$}}
\put(33,16){\circle*{1}}\put(27,28){\circle*{1}}
\filltype{black}
\put(31,20){\circle*{1}}\put(32,20){\makebox(0,0)[l]{\small$s_i$}}
\put(28,26){\circle*{1}}\put(29,26){\makebox(0,0)[l]{\small$s_j$}}
\put(30,34){\circle*{1}}\put(31,34){\makebox(0,0)[l]{\small$s_k$}}
\put(30,10){\circle*{1}}\put(31,10){\makebox(0,0)[l]{\small$s_l$}}
\put(50,22){\vector(1,0){20}}\put(60,23){\makebox(0,0)[b]{$p_i$}}
\put(140,22){\makebox(0,0)[r]{\large$\P^1_y$}}
\drawline(110,6)(110,38)
\filltype{white}
\put(110,10){\circle*{1}}
\put(108,10){\makebox(0,0)[r]{\small$p_i\circ s_-=(1:0)$}}
\put(110,34){\circle*{1}}
\put(108,34){\makebox(0,0)[r]{\small$p_i\circ s_+\!=\!(0:1)$}}
\filltype{black}
\put(110,22){\circle*{1}}
\put(108,22){\makebox(0,0)[r]{\small$p_i\circ s_i\!=\!(1:1)$}}
\put(110,28){\circle*{1}}
\put(112,28){\makebox(0,0)[l]{\small$p_i\circ s_j$}}
\put(112,34){\makebox(0,0)[l]{\small$p_i\circ s_k\!=\!p_i\circ s_+$}}
\put(112,10){\makebox(0,0)[l]{\small$p_i\circ s_l\!=\!p_i\circ s_-$}}
\end{picture}
\vspace{-11mm}
\caption{}
\vspace{-1mm}
\end{figure}

(2) {\it Contraction with respect to two sections onto an $A_1$-curve.}
Let $s_{i_1},s_{i_2}$ be two of the sections of $C$. Then there is 
a contraction morphism $p_{\{i_1,i_2\}}\colon C\to C_{\{i_1,i_2\}}$
and $C_{\{i_1,i_2\}}$ is an $A_1$-curve over $Y$.
The curve $C_{\{i_1,i_2\}}$ contains the information about the
relative positions of the sections $s_{i_1},s_{i_2}$ in $C$. 
This data relates the two contraction morphisms with respect to 
the sections $s_{i_1}$ and $s_{i_2}$.

\begin{con}\label{con:A_1curve-A_1data}
Let $(C\to Y,s_-,s_+,s_1,s_2)$ be an $A_1$-curve. Then we 
have the morphism $p_2\colon C\to\P^1_Y$ for the section 
$s_2$ such that, with respect to homogeneous coordinates 
$z_0^{(2)},z_1^{(2)}$, the sections $s_-,s_+,s_2$ become the 
$(1:0),(0:1),(1:1)$-section of $\P^1_Y$. The section $s_1$ 
determines $(t_{1,2}:t_{2,1}):=p_2\circ s_1$. This can be 
rewritten as a line bundle with two generating sections 
$(\mathscr L_{\{1,2\}},\{t_{1,2},t_{2,1}\})$ up to isomorphism.

Equivalently, we can consider the morphism $p_1\colon C\to\P^1_Y$ 
for the section $s_1$. Then $s_-,s_+,s_1$ become the 
$(1:0),(0:1),(1:1)$-section with respect to homogeneous coordinates 
$z_0^{(1)},z_1^{(1)}$ and $p_1\circ s_2=(t_{2,1}:t_{1,2})$.
\end{con}

\begin{lemma}\label{le:A_1-data}
The morphism $p_1\times p_2\colon C\to\P^1_Y\times\P^1_Y$ maps 
$C$ to the closed subvariety given by the homogeneous equation
\[t_{1,2}z_1^{(2)}z_0^{(1)}=t_{2,1}z_0^{(2)}z_1^{(1)}\,.\] 
\end{lemma}
\begin{proof}
It suffices to show the lemma for the strata of $Y$ corresponding 
to the three possible combinatorial types of $A_1$-curves.
On the strata corresponding to reducible curves the above
equation becomes $z_1^{(2)}z_0^{(1)}=0$ (resp.\
$z_0^{(2)}z_1^{(1)}=0$). These equations are satisfied, 
e.g.\ in the first case $z_1^{(2)}=0$ on the component 
containing $s_-,s_1$ and $z_0^{(1)}=0$ on the component 
containing $s_+,s_2$. 
On the remaining stratum $Y'\subseteq Y$ both $(z_1^{(1)}:z_0^{(1)})$ 
and $(z_1^{(1)}:z_0^{(1)})$ are homogeneous coordinates of $\P^1_{Y'}$
and related by $(t_{1,2}z_1^{(2)}:t_{2,1}z_0^{(2)})=(z_1^{(1)}:z_0^{(1)})$.
\end{proof}

\smallskip

(3) {\it Contraction with respect to three sections onto an $A_2$-curve.}
Let $s_{i_1},s_{i_2},s_{i_3}$ be three of the sections of $C$. Then 
there is a contraction morphism $p_{\{i_1,i_2,i_3\}}\colon C\to 
C_{\{i_1,i_2,i_3\}}$ and $C_{\{i_1,i_2,i_3\}}$ is an $A_2$-curve 
over $Y$. This curve contains the information about the relative 
positions of the pairs of two sections in a set of three sections 
$s_{i_1},s_{i_2},s_{i_3}$ of $C$. These data are related by one equation:

\begin{lemma}\label{le:A_2-data}
Let $(C,s_-,s_+,s_1,s_2,s_3)$ be an $A_2$-curve over $Y$. 
Then the collection of $A_1$-data 
$\{(\mathscr L_{\{1,2\}},\{t_{1,2},t_{2,1}\}),
(\mathscr L_{\{2,3\}},\{t_{2,3},t_{3,2}\}),
(\mathscr L_{\{3,1\}},\{t_{3,1},t_{1,3}\})\}$
extracted by the method of Construction {\rm\ref{con:A_1curve-A_1data}} 
from the $A_1$-curves $C_{\{i_1,i_2\}},C_{\{i_2,i_3\}},C_{\{i_3,i_1\}}$ 
is $A_2$-data, i.e.\ the sections satisfy the equation 
\[t_{1,2}t_{2,3}t_{3,1}=t_{2,1}t_{3,2}t_{1,3}\]
in $\mathscr L_{\{1,2\}}\otimes\mathscr L_{\{2,3\}}\otimes
\mathscr L_{\{3,1\}}$.
\end{lemma}
\begin{proof}
Again, this equation is satisfied on a closed subset and it suffices 
to consider the situation for the strata corresponding to the 
different combinatorial types.
If the sections are not contained in the same component, then
both sides vanish.
Over the remaining stratum $Y'\subseteq Y$ we have a bundle $\P^1_{Y'}$.
Using the formula  
$(t_{i,j}z_1^{(j)}:t_{j,i}z_0^{(j)})=(z_1^{(i)}:z_0^{(i)})$
for the homogeneous coordinates we obtain
$(z_1^{(1)}:z_0^{(1)})=(t_{1,2}z_1^{(2)}:t_{2,1}z_0^{(2)})
=(t_{1,2}t_{2,3}z_1^{(3)}:t_{2,1}t_{3,2}z_0^{(3)})
=(t_{1,2}t_{2,3}t_{3,1}z_1^{(1)}:t_{2,1}t_{3,2}t_{1,3}z_0^{(1)})$,
so $(t_{1,2}t_{2,3}t_{3,1}:t_{2,1}t_{3,2}t_{1,3})=(1:1)$.
\end{proof}

\smallskip

For an $A_n$-curve $(C\to Y,s_-,s_+,s_1,\ldots,s_n)$ we have
the contraction morphisms $p_i\colon C\to\P_Y(\pi_*\mathscr\O_C(s_i))
\cong(\P^1_{\{\pm\alpha_i\}})_Y$ where on $(\P^1_{\{\pm\alpha_i\}})_Y$
we have homogeneous coordinates $z_{-\alpha_i},z_{\alpha_i}$ 
such that, in these coordinates, $s_-,s_+,s_i$ become the 
$(1:0),(0:1),(1:1)$-section of $\P^1_Y$. We have the roots
$\alpha_i=u_i-u_{n+2}$ and $\beta_{ij}=\alpha_i-\alpha_j=u_i-u_j$.

\smallskip

\begin{thm}
There is an isomorphism between the functor $F_{A_n}$ and the 
moduli functor of $A_n$-curves $\overline{L}_{n+1}$
such that the universal $A_n$-data on $X(A_n)$ is mapped to 
the universal $A_n$-curve over $X(A_n)$.
\end{thm}
\begin{proof}
Let $Y$ be a scheme. For $A_n$-data on $Y$ we construct an 
$A_n$-curve $C$ over $Y$ via equations in $P(A_{n+1}/A_n)_Y=
\prod_{i=1}^{n+1}(\P^1_{\{\pm\alpha_i\}})_Y$ as in Remark 
\ref{rem:emb-univcurve} with the given $A_n$-data on $Y$ 
replacing the universal $A_n$-data on $X(A_n)$.
This is an $A_n$-curve since all $A_n$-data are pull-back of the
universal $A_n$-data on $X(A_n)$ and so the constructed curve
is a pull-back of the universal $A_n$-curve over $X(A_n)$.

In the other direction, given an $A_n$-curve on $Y$ we
extract $A_n$-data. For each pair of distinct sections 
$s_i,s_j$ we have a contraction morphism 
$C\to C_{\{i,j\}}$ onto an $A_1$-curve over $Y$.
From $(C_{\{i,j\}},s_i,s_j)$ we extract $A_1$-data 
$(\mathscr L_{\{\pm\beta_{ij}\}},t_{\beta_{ij}},
t_{-\beta_{ij}}):=(\mathscr L_{\{i,j\}},t_{i,j},t_{j,i})$ 
as in Construction \ref{con:A_1curve-A_1data}.
That the collection of all these data forms $A_n$-data
follows from Lemma \ref{le:A_2-data}.

\pagebreak
Both constructions commute with base-change and thus define 
morphisms of functors $F_{A_n}\to\overline{L}_{n+1}$
and $\overline{L}_{n+1}\to F_{A_n}$. 
We show that they are inverse to each other.

\smallskip

Starting with an $A_n$-curve $C$ over $Y$, we extract data
$(\mathscr L_{\{\pm\beta_{ij}\}},\{t_{-\beta_{ij}},
t_{\beta_{ij}}\})_{\{\pm\beta_{ij}\}\subseteq A_n}$ and from these 
$A_n$-data we construct an $A_n$-curve $C'\subseteq P(A_{n+1}/A_n)_Y$
as in Remark \ref{rem:emb-univcurve}.
We show that the product of the contraction morphisms $p_i$ 
defines an isomorphism of $A_n$-curves $C\to C'\subseteq P(A_{n+1}/A_n)_Y$.
The morphism $\prod_ip_i\colon C\to P(A_{n+1}/A_n)_Y$ factors through 
the inclusion $C'\subseteq P(A_{n+1}/A_n)_Y$ by Lemma \ref{le:A_1-data}. 
The morphism $C\to C'$ is an isomorphism, because it is an isomorphism 
on the fibres by Proposition \ref{prop:emb-A_n-curve} and the curves 
are flat over $Y$.
The sections and the involution of $C'$ were defined such that 
$C\to C'$ is an isomorphism of $A_n$-curves.

\smallskip

Starting with $A_n$-data we construct an $A_n$-curve and extract
$A_n$-data from it. It is easily verified that this $A_n$-data 
coincides with the original $A_n$-data, since the contraction 
morphism with respect to a section of an embedded $A_n$-curve
$C'\subseteq P(A_{n+1}/A_n)_Y$ is induced by projection
onto the corresponding factor.
\end{proof}

\end{document}